\documentclass{article}
\usepackage{amsmath,amssymb}
\usepackage[all]{xy}
\usepackage{caption,graphicx}
\font\ibf=cmbxti10

\title{Equivariant, locally finite inverse representations with uniformly bounded zipping length, for arbitrary finitely presented groups}
\author{Valentin {\sc Po\'enaru}\footnote{Professor Emeritus at the Universit\'e Paris Sud-Orsay, Math\'ematiques 425, 91405 Orsay Cedex, France. e-mail: valpoe@hotmail.com}}
\date{February 2009}

\begin{document}

\maketitle

\vglue 1cm

\section{Introduction}\label{sec0}
\setcounter{equation}{0}

This is the first of a series of papers the aim of which is to give complete proofs for the following

\bigskip

\noindent {\bf Theorem.} {\it Any finitely presented group $\Gamma$ is QSF.}

\bigskip

This has already been announced in my preprint \cite{24}, where a few more details about the general context of the theorem above are also given. The QSF ($=$ quasi simply filtered) for locally compact spaces $X$ and/or for finitely presented groups $\Gamma$, has been previously introduced by S.~Brick, M.~Mihalik and J.~Stallings \cite{3}, \cite{29}.

\smallskip

A locally compact space $X$ is said to be QSF, and it is understood that from now on we move in the simplicial category, if the following thing happens. For any compact $k \subset X$ we can find an abstract compact $K$, where by ``abstract'' we mean that $K$ is not just another subspace of $X$, such that $\pi_1 K = 0$, and coming with an embedding $k \overset{j}{\longrightarrow} K$, such that there is a continuous map $f$ entering the following commutative diagram
$$
\xymatrix{
k \ar[rr]^{j} \ar@{_{(}->}[dr] &&K \ar[dl]^{f}  \\ 
&X
}
$$
with the property that $M_2 (f) \cap jk = \emptyset$. Here $M_2 (f) \subset K$ is the set of double points $x \in K$ s.t. ${\rm card} \, f^{-1} f(x) > 1$. [We will also use the notation $M^2 (f) \subset K \times K$ for the pairs $(x,y)$, $x \ne y$, $fx = fy$.] One of the virtues of QSF is that if $\Gamma$ is a finitely presented group and $P$ is a PRESENTATION for $\Gamma$, {\it i.e.} a finite complex with $\pi_1 P = \Gamma$, then ``$\tilde P \in \Gamma$'' is presentation independent, {\it i.e.} if one presentation of $\Gamma$ has this property, then all presentations of $\Gamma$ have it too. Also, if this is the case, then we will say that $\Gamma$ itself is QSF. The QSF {\ibf is} a bona fide group theoretical concept. Another virtue of QSF is that if $V^3$ is an open contractible 3-manifold which is QSF, then $\pi_1^{\infty} V^3 = 0$.

\smallskip

The special case when $\Gamma = \pi_1 M^3$, with $M^3$ a closed 3-manifold, of the theorem above, is already known indeed, since that particular case is already a consequence of the full big work of G.~Perelman on the Thurston Geometrization of 3-manifolds \cite{11}, \cite{12}, \cite{13}, \cite{9}, \cite{2}, \cite{1}, \cite{8}. This is one of the impacts of the work in question on geometric group theory. Now, for technical reasons which will become clear later on, when flesh and bone will be put onto $\Gamma$ by picking up a specific presentation $P$, then we will make the following kind of choice for our $P$. We will disguise our general $\Gamma$ as a $3$-manifold group, in the sense that we will work with a compact bounded {\ibf singular} $3$-manifold $M^3 (\Gamma)$, with $\pi_1 M^3 (\Gamma) = \Gamma$. The $M^3 (\Gamma)$'s will be described with more detail in the section~\ref{sec1} below. Keep in mind that $M^3 (\Gamma)$ is NOT a smooth manifold and it certainly has no vocation of being geometrized in any sense.

\smallskip

A central notion for our whole approach is that of REPRESENTATION for $\Gamma$ and/or for $\tilde M^3 (\Gamma)$, two objects which up to quasi-isometry are the same. Our terminology here is quite unconventional, our representations proceed like ``$\longrightarrow \Gamma$'', unlike the usual representation of groups which proceed like ``$\Gamma \longrightarrow$''. The reason for this slip of tongue will be explained later. Here is how our representations are defined. We have to start with some infinite, not necessarily locally finite complex $X$ of dimension $2$ or $3$ (but only $\dim X = 3$ will occur in the present paper), endowed with a non-degenerate cellular map $X \overset{f}{\longrightarrow} \tilde M^3 (\Gamma)$.

\smallskip

Non-degenerate means here that $f$ injects on the individual cells. The points of $X$ where $f$ is not immersive, are by definition the {\ibf mortal singularities} ${\rm Sing} \, (f) \subset X$. We reserve the adjective {\ibf immortal} for the singularities of $\tilde M^3 (\Gamma)$ and/or of $M^3 (\Gamma)$.

\smallskip

Two kinds of equivalence relations will be associated to such non-degenerate maps $f$, namely the
$$
\Psi (f) \subset \Phi (f) \subset X \times X \, .
$$
Here $\Phi (f)$ is the simple-minded equivalence relation where $(x,y) \in \Phi (f)$ means that $f(x) = f(y)$. The more subtle and not so easily definable $\Psi (f)$, is the smallest equivalence relation compatible with $f$, which kills all the mortal singularities. More details concerning $\Psi (f)$ will be given in the first section of this present paper. But before we can define our representations for $\Gamma$ we also have to review the notion of GSC, {\it i.e.} {\ibf ``geometrically simply connected''}. This stems from differential topology where a smooth manifold is said to be GSC iff it has a handle-body decomposition without handles of index $\lambda = 1$, and/or such that the $1$-handles and $2$-handles are in cancelling position, and a certain amount of care is necessary here, since what we are after are rather non-compact manifolds with non-empty boundary. As it is explained in \cite{24} there are some deep connections between the concept GSC from differential topology (for which, see also \cite{20}, \cite{21}, \cite{25}) and the concept QSF in group theory.

\smallskip

But then, the GSC concept generalizes easily for cell-complexes and the ones of interest here will always be infinite. We will say that our $X$ which might be a smooth non-compact manifold (possibly with boundary $\ne \emptyset$) or a cell-complex, like in $X \longrightarrow \tilde M^3 (\Gamma)$, is GSC if it has a cell-decomposition or handle-body decomposition, according to the case of the following type. Start with a PROPERLY embedded tree (or with a smooth regular neighbourhood of such a tree); with this should come now a cell-decomposition, or handle-decomposition
$$
X = T \cup \sum_{1}^{\infty} \left\{ \mbox{$1$-handles (or $1$-cells)} \, H_i^1 \right\} \cup \sum_{1}^{\infty} H_j^2 \cup \left\{ \sum_{k;\lambda \geq 2} H_k^{\lambda} \right\} \, ,
$$
where the $i$'s and the $j$'s belong to a same countable set, such that the {\ibf geometric intersection matrix}, which counts without any kind of $\pm$ signs, how many times $H_j^2$ goes through $H_i^1$, takes the following ``easy'' id $+$ nilpotent form
\begin{equation}
\label{eq0.1}
H_j^2 \cdot H_i^1 = \delta_{ji} + a_{ji} \, , \quad \mbox{where} \quad a_{ji} \in Z_+ \quad \mbox{and} \quad a_{ji} > 0 \Rightarrow j > i \, .
\end{equation}

In this paper we will distinguish between PROPER, meaning inverse image of compact is compact, and proper meaning (inverse image of boundary) $=$ boundary.

\smallskip 

Also, there is a notion of ``difficult'' id $+$ nilpotent, gotten by reversing the last inequality in (\ref{eq0.1}) and this is no longer GSC. For instance, the Whitehead manifold ${\rm Wh}^3$, several times mentioned in this sequence of papers and always impersonating one of the villains of the cast, admits handle-body decompositions of the difficult id $+$ nilpotent type. But then, ${\rm Wh}^3$ is certainly not GSC either. Neither are the various manifolds ${\rm Wh}^3 \times B^n$, $n \geq 1$ (see here \cite{15}).

\smallskip

We can, finally, define the REPRESENTATIONS for $\tilde M^3 (\Gamma)$ (and/or $\Gamma$). These are non-degenerate maps, like above
\begin{equation}
\label{eq0.2}
X \overset{f}{\longrightarrow} \tilde M^3 (\Gamma)
\end{equation}
with the following features

\medskip

I) $X \in {\rm GSC}$; we call it the representation space

\medskip

II) $\Psi (f) = \Phi (f)$

\medskip

III) $f$ is {\ibf ``essentially surjective''} which, for the case when $\dim X = 3$ this simply means that $\overline{{\rm Im} f} = \tilde M^3 (\Gamma)$.

\medskip

Here are some comments: A) The fact that it is a {\ibf group} $\Gamma$ which is being ``represented'' this way, does {\ibf not} occur explicitly in the features I), II), III). The point is that, when groups $\Gamma$ ({\it i.e.} $\tilde M^3 (\Gamma)$) are being represented, this opens the possibility for the representation to be {\ibf equivariant}, a highly interesting proposition, as it will turn out soon.

\smallskip

On the other hand, as it will soon become clear, whenever a map $(X,f)$ has the features I), II), III), then this {\it forces} whatever sits at the target of $f$ to be simply-connected. So, a priori at least, one may try to represent various simply-connected objects. Long ago, the present author has started by representing this way homotopy $3$-spheres $\Sigma^3$ (see here \cite{16} and \cite{6}, papers which eventually culminated with \cite{23}). Then, I also represented universal covering spaces of smooth closed $3$-manifolds (\cite{18}, \cite{19}) or even the wild Whitehead manifold (see \cite{27}). My excuse for calling the thing occuring in (\ref{eq0.2}) a ``representation'' for $\Gamma$, is that I had started to use this terminology already long ago, in contexts where no group was present; and by now the sin has been committed already.

\smallskip

Initially, the present kind of representations were called ``pseudo-spine representations'' (see \cite{16}, for instance). But today, in contexts where no confusion is possible I will just call them ``representations''. Only once, in the title to the present paper, I have also added the adjective ``inverse'', with the sole purpose of avoiding confusion.

\medskip

\noindent B) We will call $n = \dim X$ the {\ibf dimension} of the representation, which can be $n=3$, like it will be the case now, or $n=2$ which is the really useful thing, to be developed in the next papers. The reason for using the presentation $M^3 (\Gamma)$ is to be able to take advantage of the richness of structure of $M_2 (f)$ when maps $(\dim = 2) \overset{f}{\longrightarrow} (\dim = 3)$ are concerned.

\medskip

\noindent C) When $\dim X = 2$, then the ``essential surjectivity'' from III) will mean that $\tilde M^3 (\Gamma) - \overline{{\rm Im} (f)}$ is a countable union of small, open, insignificant ($=$ cell-like) subsets. This ends our comments.

\medskip

Let us come back now to the condition $\Psi (f) = \Phi (f)$ from II) above and give it a more geometric meaning. We consider then
$$
\Phi (f) \supset M^2 (f) \subset X \times X \supset {\rm Sing} \, (f)
$$
where, strictly speaking, ``${\rm Sing} \, (f)$'' should be read ``${\rm Diag} \, ({\rm Sing} \, (f))$''. We will extend $M^2 (f)$ to $\hat M^2 (f) = M^2 (f) \cup {\rm Sing} \, (f) \subset X \times X$. With this, the condition $\Psi (f) = \Phi (f)$ means that, at the level of $\hat M^2 (f)$, any $(x,y) \in M^2 (f)$ can be joined by continuous paths to the singularities. The figure~1.1 below should be enough, for the time being, to make clear what we talk about now; more formal definitions will be given later. Anyway, the kind of thing which figure~1.1 displays will be called a zipping path (or a zipping strategy) $\lambda (x,y)$ for $(x,y) \in M^2 (f)$. Here is a way to make sense of the length $\Vert \lambda (x,y) \Vert$ of such a zipping path.

$$
\includegraphics[width=11cm]{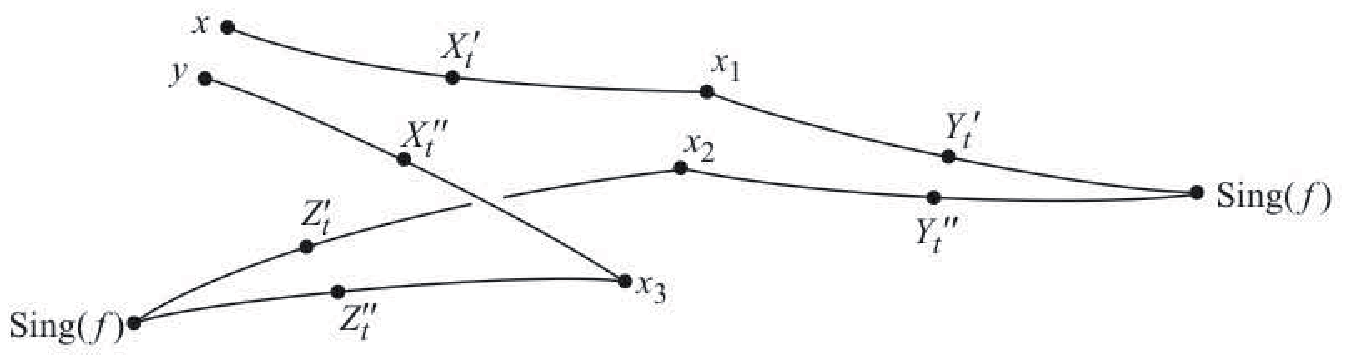}
$$
\centerline{Figure 1.1. A zipping path $\lambda (x,y)$ for $(x,y) \in M^2 (f)$.}
\begin{quote}
Here $(x_1 , x_2 , x_3) \in M^3 (f)$ ($=$ triple points) and the various moving points $(X'_t , X''_t)$, $(Y'_t , Y''_t)$, $(Z'_t , Z''_t)$ which depend continuously on $t \in [0,1]$, are in $M^2 (f)$. The condition $\Psi (f) = \Phi (f)$ means the existence of a zipping path, like above, for any double point. But, of course, zipping paths are {\ibf not} unique.
\end{quote}

Although $M^3 (\Gamma)$ is not a smooth manifold, riemannian metrics can be defined for it. Such metrics can then be lifted to $\tilde M^3 (\Gamma)$, to $X$ or to $X \times X$. Using them we can make sense of $\Vert \lambda (x,y) \Vert$ which, of course, is only well-defined up to quasi-isometry. The condition $\Psi (f) = \Phi (f)$ also means that the map $f$ in (\ref{eq0.2}) is realizable (again not uniquely, of course), via a sequence of folding maps. Again, like in the context of figure~1.1, this defines a zipping or a zipping strategy.

\smallskip

When a double point $(x,y) \in M^2 (f)$ is given, we will be interested in quantities like ${\rm inf} \, \Vert \lambda (x,y) \Vert$, taken over all the possible continuous zipping paths for $(x,y)$. A quasi-isometrically equivalent, discrete definition, would be here to count the minimum member of necessary folding maps for closing $(x,y)$. But we will prefer the continuous version, rather than the discrete one. $\Box$

\bigskip

The object of the present paper is to give a proof for the following

\bigskip

\noindent REPRESENTATION THEOREM. {\it For any $\Gamma$ there is a $3$-dimensional representation
\begin{equation}
\label{eq0.3}
\xymatrix{
Y(\infty) \ar[rr]^{\!\!\!\!\!\!\!\!\!\!\ g(\infty)}  && \ \tilde M^3 (\Gamma) \, , 
} 
\end{equation}
such that the following things should happen too}

\medskip

1) {\it The $Y(\infty)$ is {\ibf locally finite.}}

\medskip

2) {\it The representation} (\ref{eq0.3}) {\it is EQUIVARIANT. Specifically, the representation space $Y(\infty)$ itself comes equipped with a free action of $\Gamma$ and then, for each $x \in Y(\infty)$ and $\gamma \in \Gamma$, we have also that}
$$
g(\infty) (\gamma x) = \gamma g (\infty)(x) \, .
$$

3) {\it There is a {\ibf uniform bound} $M > 0$ such that, for any $(x,y) \in M^2 (g(\infty))$ we have}
\begin{equation}
\label{eq0.4}
\underset{\lambda}{\rm inf} \, \Vert \lambda (x,y) \Vert < M \, ,
\end{equation}
{\it when $\lambda$ runs over all zipping paths for $(x,y)$.}

\bigskip

Here are some comments concerning our theorem above.

\medskip

\noindent D) In general, representation spaces are violently {\ibf non} locally finite. So, the feature 1) is pretty non-trivial. In our present case of a $3$-dimensional representation, the image $g(\infty) \, Y(\infty)$ will automatically be the interior of the singular manifold $\tilde M^3 (\Gamma)$ which, of course, is locally finite. This does not involve in any way the finer aspects of the representation theorem, namely locally finite source, equivariance and bounded zipping length. But then, also, as we shall see below, when we go to $2^{\rm d}$ representation, even with our three finer features present, the image is no longer locally finite, generally speaking at least. This is a major issue for this sequence of papers.

\medskip

\noindent E) It is in 2) above that the group property is finally used.

\medskip

\noindent F) In the special case when $\Gamma = \pi_1 M^3$ for a smooth closed $3$-manifold, parts 1) + 2) of the REPRESENTATION THEOREM are already proved in \cite{26}. It turns out that the way the paper \cite{26} was constructed, makes it very easy for us to extend things from $M_{\rm smooth}^3$ to the present $M^3 (\Gamma)_{\rm SINGULAR}$. End of comments.

\medskip

The first section of this paper reviews the $\Psi/\Phi$ theory in our present context, while the section II $+$ III give the proof of the REPRESENTATION THEOREM above. They rely heavily on \cite{26}. In the rest of the present introduction, we will start by giving a rather lengthy sketch of how the REPRESENTATION THEOREM fits into the proof of the main theorem stated in the very beginning, {\it i.e.} the fact that all finitely presented $\Gamma$'s are QSF. This  will largely overlap with big pieces of the announcement \cite{24}. We will end the introduction by briefly reviewing some possible extensions of the present work.

\bigskip

\noindent AN OVERVIEW OF THE PROOF THAT ``ALL $\Gamma$'s ARE QSF''. The representation theorem above provides us with a $3^{\rm d}$ representation space $Y(\infty)$. But we actually need a $2$-dimensional representation, in order to be able to go on. So we will take a sufficiently dense skeleton of the $Y(\infty)$ from the representation theorem and, making use of it, we get the next representation theorem below.

\bigskip

\noindent $2$-DIMENSIONAL REPRESENTATION LEMMA. {\it For any finitely presented $\Gamma$ there is a $2$-dimensional representation
\begin{equation}
\label{eq0.5}
X^2 \overset{f}{\longrightarrow} \tilde M^3 (\Gamma)
\end{equation}
which, just like in the preceeding $3^{\rm d}$ REPRESENTATION THEOREM, is}

\medskip

1) {\it Locally finite; we will call this the first finiteness condition.}

\medskip

2) {\it Equivariant;}

\medskip

3) {\it With uniformly bounded zipping length and, moreover, such that we also have the next features below.}

\medskip

4) {\it (The second finiteness condition.) For any tight compact transversal $\Lambda$ to $M_2 (f) \subset X^2$ we have}
\begin{equation}
\label{eq0.6}
\mbox{card} \, (\lim \, (\Lambda \cap M_2 (f))) < \infty \, . 
\end{equation}

\medskip

5) {\it The closed subset}
\begin{equation}
\label{eq0.7}
\mbox{LIM} \, M_2 (f) \underset{\rm def}{=} \ \bigcup_{\Lambda} \ \lim \, (\Lambda \cap M_2 (f)) \subset X^2 
\end{equation}
{\it is a locally finite graph and $f \, {\rm LIM} \, M_2 (f) \subset f X^2$ is also a {\ibf closed} subset.}

\medskip

6) {\it Let $\Lambda^*$ run over all tight transversals to ${\rm LIM} \, M_2 (f)$. Then we have}
$$
\bigcup_{\Lambda^*} \ (\Lambda^* \cap M_2 (f)) = M_2 (f) \, .
$$

\bigskip

One should notice that ${\rm LIM} \, M_2 (f) = \emptyset$ is equivalent to $M_2 (f) \subset X^2$ being a closed subset. It so happens that, if this is the case, then it is relatively easy to prove that $\Gamma \in $ QSF. This makes that in all this sequence of papers we will only consider the situation ${\rm LIM} \, M_2 (f) \ne \emptyset$, {\it i.e.} the worst possible case.

\smallskip

Now, once $M_2 (f)$ is NOT a closed subset, the (\ref{eq0.6}) is clearly the next best option. But then, in \cite{27} it is shown that if, forgetting about groups and group actions, we play this same game, in the most economical manner for the Whitehead manifold ${\rm Wh}^3$, instead of $\tilde M^3 (\Gamma)$, then generically the set $\lim (\Lambda \cap M_2 (f))$ becomes a Cantor set, which comes naturally equipped with a feedback loop, turning out to be directly related to the one which generates the Julia sets in the complex dynamics of quadratic maps. See \cite{27} for the exact Julia sets which are concerned here.

\smallskip

Continuing with our list of comments, in the context of 6) in our lemma above, there is a transversal holonomy for ${\rm LIM} \, M_2 (f)$, which is quite non trivial. Life would be easier without that.

\smallskip

The $X^2$ has only mortal singularities, while the $fX^2$ only has immortal ones. [In terms of \cite{16}, \cite{6}, the singularities of $X^2$ (actually of $f$) are all ``undrawable singularities''. At the source, the same is true for the immortal ones too.]

\smallskip

The representation space $X^2$ in Lemma~5 is locally finite but, as soon as ${\rm LIM} \, M_2 (f) \ne \emptyset$ (which, as the reader will retain, is the main source of headaches in the present series of papers), the $fX^2 \subset \tilde M^3 (\Gamma)$ is {\ibf not}.

\smallskip

At this point, we would like to introduce canonical smooth high-dimensional thickenings for this $fX^2$. Here, even if we temporarily forget about $f \, {\rm LIM} \, M_2 (f)$, there are still the immortal singularities (the only kind which $f X^2$ possesses), and these accumulate at finite distance. This is something which certainly has to be dealt with but, at the level of this present sketch we will simply ignore it. Anyway, even dealing with isolated immortal singularities only, requires a relatively indirect procedure. One starts with a 4-dimensional smooth thickening (remember that, provisionally, we make as if $f \, {\rm LIM} \, M_2 (f) = \emptyset$). This 4$^{\rm d}$ thickening is not uniquely defined, it depends on a desingularization (see \cite{6}). Next one takes the product with $B^m$, $m$ large, and this washes out the desingularization-dependence.

\medskip

In order to deal with ${\rm LIM}Ê\, M_2 (f) \ne \emptyset$, the procedure sketched above has to be supplemented with appropriate {\ibf punctures}, by which we mean here pieces of ``boundary'' of the prospective thickened object which are removed, or sent to infinity. This way we get a more or less canonical {\ibf smooth} high-dimensional thickening for $fX^2$, which we will call $S_u (\tilde M^3 (\Gamma))$. Retain that the definition of $S_u (\tilde M^3 (\Gamma))$ has to include at least a first batch of punctures, just so as to get a smooth object. By now we can also state the following

\bigskip

\noindent {\bf Main Lemma.} {\it $S_u (\tilde M^3 (\Gamma))$ is} GSC.

\bigskip

We will go to some length with the sketch of proof for this Main Lemma below, but here are some comments first. What our main lemma says, is that a certain smooth, very high-dimensional non-compact manifold with large non-empty boundary, very much connected with $\tilde M^3 (\Gamma)$ but without being exactly a high-dimensional thickening of it, is GSC. We will have a lot to say later concerning this manifold, our $S_u (\tilde M^3 (\Gamma))$. For right now, it suffices to say that it comes naturally equipped with a free action of $\Gamma$ but for which, unfortunately so to say, the fundamental domain $S_u (\tilde M^3 (\Gamma))/\Gamma$ fails to be compact. If one thinks of the $\tilde M^3 (\Gamma)$ itself as being made out of fundamental domains which are like solid compact ice-cubes, then the $S_u (\tilde M^3 (\Gamma))$ is gotten, essentially, by replacing each of these compact ice-cubes by a non-compact infinitely foamy structure, and then going to a very high-dimensional thickening of this foam.

\smallskip

In order to clinch the proof of ``$\forall \, \Gamma \in {\rm QSF}$'' one finally also needs to prove the following kind of implication

\bigskip

\noindent {\bf Second Lemma.} {\it We have the following implication}
$$
\{ S_u (\tilde M^3 (\Gamma)) \in {\rm GSC}, \ \mbox{which is what the MAIN LEMMA claims}\} \Longrightarrow 
$$
\begin{equation}
\label{eq0.8}
\Longrightarrow \{ \Gamma \in {\rm QSF}\} \, .
\end{equation}

\bigskip

At the level of this short outline there is no room even for a sketch of proof of the implication (\ref{eq0.8}) above. Suffices to say here that this (\ref{eq0.8}) is a considerably easier step than the main lemma itself, about which a lot will be said below.

\smallskip

We end this OVERVIEW with a SKETCH OF THE PROOF OF THE MAIN LEMMA.

\smallskip

For the lemma to be true a second batch of punctures will be necessary. But then, we will also want our $S_u (\tilde M^3 (\Gamma))$ to be such that we should have the implication (\ref{eq0.8}). This will turn out to put a very strict {\ibf ``Stallings barrier''} on how much punctures the definition of $S_u (\tilde M^3 (\Gamma))$ can include, at all. 

\smallskip

So, let $Y$ be some low-dimensional object, like $fX^2$ or some 3-dimensional thickening of it. When it is a matter of punctures, these will be put into effect directly at the level of $Y$ (most usually by an appropriate infinite sequence of Whitehead dilatations), before going to higher dimensions, so that the high dimensional thickening should be {\ibf transversally compact} with respect to the $Y$ to be thickened. The reason for this requirement is that when we will want to prove (\ref{eq0.8}), then arguments like in \cite{15} will be used (among others), and these ask for transversal compactness.  For instance, if $V$ is a low-dimensional non-compact manifold, then
$$
\mbox{$V \times B^m$ is transversally compact, while $V \times \{ B^m$ with some boundary}
$$
$$
\mbox{punctures$\}$, or even worst $V \times R^m$, is not.}
$$

This is our so-called ``Stallings barrier'', putting a limit to how much punctures we are allowed to use. The name is referring here to a corollary of the celebrated Engulfing Theorem of John Stallings, saying that with appropriate dimensions, if $V$ is an {\ibf open} contractible manifold, then $V \times R^p$ is a Euclidean space. There are, of course, also infinitely more simple-minded facts which give GSC when multiplying with $R^p$ or $B^p$.

\smallskip

We take now a closer look at the $n$-dimensional smooth manifold $S_u (\tilde M^3 (\Gamma))$ occuring in the Main Lemma. Here $n=m+4$ with high $m$. Essentially, we get our $S_u (\tilde M^3 (\Gamma))$ starting from an initial GSC low-dimensional object, like $X^2$, by performing first a gigantic quotient-space operation, namely our zipping, and finally thickening things into something of dimension $n$. But the idea which we will explore now, is to construct another smooth $n$-dimensional manifold, essentially starting from an already $n$-dimensional GSC smooth manifold, and then use this time a gigantic collection of additions and inclusion maps, rather than quotient-space projections, which should somehow mimic the zipping process. We will refer to this kind of thing as a {\ibf geometric realization of the zipping}. The additions which are allowed in this kind of realization are Whitehead dilatations, or adding handles of index $\lambda \geq 2$. The final product of the geometric realization, will be another manifold of the same dimension $n$ as $S_u (\tilde M^3 (\Gamma))$, which we will call $S_b (\tilde M^3 (\Gamma))$ and which, {\it a priori} could be quite different from the $S_u (\tilde M^3 (\Gamma))$. To be more precise about this, in a world with ${\rm LIM} \, M_2 (f) = \emptyset$ we would quite trivially have $S_u = S_b$ but, in our real world with ${\rm LIM} \, M_2 (f) \ne \emptyset$, there is absolutely no {\it a priori} reason why this should be so. Incidentally, the subscripts ``$u$'' and ``$b$'', refer respectively to ``usual'' and ``bizarre''.

\smallskip

Of course, we will want to compare $S_b (\tilde M^3 (\Gamma))$ and $S_u (\tilde M^3 (\Gamma))$. In order to give an idea of what is at stake here, we will look into the simplest possible local situation with ${\rm LIM} \, M_2 (f)$ present. Ignoring now the immortal singularities for the sake of the exposition, we consider a small smooth chart $U = R^3 = (x,y,z)$ of $\tilde M^3 (\Gamma)$, inside which live $\infty + 1$ planes, namely $W = (z=0)$ and the $V_n = (x=x_n)$, where $x_1 < x_2 < x_3 < \ldots$ with $\lim x_n = x_{\infty}$. Our local model for $X^2 \overset{f}{\longrightarrow} \tilde M^3 (\Gamma)$ is here $f^{-1} \, U = W + \underset{1}{\overset{\infty}{\sum}} \ V_n \subset X^2$ with $f \mid \Bigl( W + \underset{1}{\overset{\infty}{\sum}} \ V_n \Bigl)$ being the obvious map. We find here that the line $(x = x_{\infty}, z=0) \subset W$ is in ${\rm LIM} \, M_2 (f)$ and the situation is sufficiently simple so that we do not need to distinguish here between ${\rm LIM} \, M_2 (f)$ and $f \, {\rm LIM} \, M_2 (f)$.

\smallskip

Next, we pick up a sequence of positive numbers converging very fast to zero $\varepsilon > \, \varepsilon_1 > \varepsilon_2 > \ldots$ and, with this, on the road to the $S_u (\tilde M^3 (\Gamma))$ from the {\ibf Main Lemma}, we will start by replacing the $f W \cup \underset{1}{\overset{\infty}{\sum}} \ V_n \subset f X^2$, with the following 3-dimensional non-compact 3-manifold with boundary
$$
M \underset{\rm def}{=} [ W \times (-\varepsilon \leq z \leq \varepsilon) - {\rm LIM} \, M_2 (f) \times \{ z = \pm \, \varepsilon \}] \ \cup
$$
\begin{equation}
\label{eq0.9}
\cup \ \sum_{1}^{\infty} V_n \times (x_n - \varepsilon_n \leq x \leq x_n + \varepsilon_n ) \, .
\end{equation}
In such a formula, notations like ``$W \times (-\varepsilon \leq z \leq \varepsilon)$'' should be read ``$W$ thickened into $-\varepsilon \leq z \leq \varepsilon$''. Here ${\rm LIM} \, M_2 (f) \times \{ \pm \, \varepsilon \}$ is a typical puncture, necessary to make our $M$ be a smooth manifold. For expository purposes, we will pretend now that $n=4$ and then $M \times [0 \leq t \leq 1]$ is a local piece of $S_u (\tilde M^3 (\Gamma))$. Now, in an ideal world (but not in ours!), the geometrical realization of the zipping process {\it via} the inclusion maps (some of which will correspond to the Whitehead dilatations which are necessary for the punctures), which are demanded by $S_b (\tilde M^3 (\Gamma))$, should be something like this. We start with the obviously GSC $n$-dimensional thickening of $X^2$, call it $\Theta^n (X^2)$; but remember that for us, here, $n=4$. Our local model should live now inside $R^4 = (x,y,z,t)$, and we will try to locate it there conveniently for the geometric realization of the zipping. We will show how we would like to achieve this for a generic section $y = $ constant.

\smallskip

For reasons to become soon clear, we will replace the normal section $y = $ const corresponding to $W$ and which should be
$$
N_y = [ - \infty < x < \infty \, , \ y = {\rm const} \, , \ - \varepsilon \leq z \leq \varepsilon \, , \ 0 \leq t \leq 1] 
$$
\begin{equation}
\label{eq0.10}
- \, (x = x_{\infty} \, , \ y = {\rm const} \, , \ z = \pm \, \varepsilon \, , \ 0 \leq t \leq 1) \, ,
\end{equation}
by the smaller $N_y \, - \, \overset{\infty}{\underset{1}{\sum}} \ {\rm DITCH} \, (n)_y$, which is defined as follows. The ${\rm DITCH} \, (n)_y$ is a thin column of height $-\varepsilon \leq z \leq \varepsilon$ and of $(x,t)$-width $4 \, \varepsilon_n$, which is concentrated around the arc
$$
(x = x_n \, , \ y = {\rm const} \, , \ - \varepsilon \leq z \leq \varepsilon \, , \ t=1) \, .
$$
This thin indentation inside $N_y$ is such that, with our fixed $y = {\rm const}$ being understood here, we should have
\begin{equation}
\label{eq0.11}
\lim_{n = \infty} \, {\rm DITCH} \, (n)_y = (x = x_{\infty} \, , \ -\varepsilon \leq z \leq \varepsilon \, , \ t=1) \, . 
\end{equation}
Notice that, in the RHS of (\ref{eq0.11}) it is exactly the $z = \pm \, \varepsilon$ which corresponds to punctures.

\medskip

Continuing to work here with a fixed, generic $y$, out of the normal $y$-slice corresponding to $V_n$, namely $(x_n - \varepsilon_n \leq x \leq x_n + \varepsilon_n \, , \ - \infty < z < \infty \, , \ 0 \leq t \leq 1)$, we will keep only a much thinner, isotopically equivalent version, namely the following
\begin{equation}
\label{eq0.12}
(x_n - \varepsilon_n \leq x \leq x_n + \varepsilon_n \, , \ - \infty < z < \infty \, , \ 1 - \varepsilon_n \leq t \leq 1) \, . 
\end{equation}
This (\ref{eq0.12}) has the virtue that it can fit now inside the corresponding ${\rm DITCH} \, (n)_y$, without touching at all the $N_y - \{{\rm DITCHES}\}$.

\smallskip

What has been carefully described here, when all $y$'s are being taken into account, is a very precise way of separating the $\infty + 1$ branches of (the thickened) (\ref{eq0.9}), at the level of $R^4$, taking full advantage of the additional dimensions ({\it i.e.} the factor $[0 \leq t \leq 1]$ in our specific case). With some work, this kind of thing can be done consistently for the whole global $fX^2$. The net result is an isotopically equivalent new model for $\Theta^n (X^2)$, which invites us to try the following naive approach for the geometric realization of the zipping. Imitating the successive folding maps of the actual zipping, fill up all the empty room left inside the ditches, by using only Whitehead dilatations and additions of handles of index $\lambda > 1$, until one has reconstructed completely the $S_u (\tilde M^3 (\Gamma))$. Formally there is no obstruction here and then also what at a single $y = $ const may look like a handle of index one, becomes ``only index $\geq 2$'', once the full global zipping is taken into account. But there {\ibf is} actually a big problem with this naive approach, {\it via} which one can certainly reconstruct $S_u (\tilde M^3 (\Gamma))$ as a set, but with the {\ibf wrong topology}, as it turns out. I will give an exact idea now of how far one can actually go, proceeding in this naive way. In \cite{28} we have, actually, tried to play the naive game to its bitter end, and the next Proposition~1.1, given here for purely pedagogical reasons, is the kind of thing one gets (and essentially nothing better than it).

\bigskip

\noindent {\bf Proposition 1.1.} {\it Let $V^3$ be any open simply-connected $3$-manifold, and let also $m \in Z_+$ be high enough. There is then an infinite collection of smooth $(m+3)$-dimensional manifolds, all of them non-compact, with very large boundary, connected by a sequence of smooth embeddings
\begin{equation}
\label{eq0.13}
X_1 \subset X_2 \subset X_3 \subset \ldots 
\end{equation}
such that}

\medskip

1) {\it $X_1$ is {\rm GSC} and each of the inclusions in {\rm (\ref{eq0.13})} is either an elementary Whitehead dilatation or the addition of a handle of index $\lambda > 1$.}

\medskip

2) {\it When one considers the union of the objects in {\rm (\ref{eq0.13})}, endowed with the {\ibf weak topology}, and there is no other, reasonable one which is usable here, call this new space $\varinjlim \, X_i$, then there is a continuous bijection}
\begin{equation}
\label{eq0.14}
\xymatrix{
\varinjlim \, X_i  \overset{\psi}{\longrightarrow} V^3 \times B^m \, . 
} 
\end{equation}

The reader is reminded that in the weak topology, a set $F \subset \varinjlim \, X_i$ is closed iff all the $F \cap X_i$ are closed. Also, the inverse of $\psi$ is not continuous here; would it be, this would certainly contradict \cite{15}, since $V^3$ may well be ${\rm Wh}^3$, for instance. This, {\it via} Brouwer, also means that $\varinjlim \, X_i$ cannot be a manifold (which would automatically be then GSC). Also, exactly for the same reasons why $\varinjlim \, X_i$ is not a manifold, it is not a metrizable space either. So, here we meet a new barrier, which I will call the {\ibf non metrizability barrier} and, when we will really realize geometrically the zipping, we better stay on the good side of it.

\smallskip

One of the many problems in this series of papers is that the Stallings barrier and the non metrizability barrier somehow play against each other, and it is instructive to see this in a very simple instance. At the root of the non metrizability are, as it turns out, things like (\ref{eq0.11}); {\it a priori} this might, conceivably, be taken care of by letting all of $(x = x_{\infty} \, , \ - \varepsilon \leq z \leq \varepsilon \, , \ t=1)$ be punctures, not just the $z = \pm \, \varepsilon$ part. But then we would also be on the wrong side of the Stallings barrier. This kind of conflict is quite typical.

\smallskip

So far, we have presented the disease and the rest of the section gives, essentially, the cure. In a nutshell here is what we will do. We start by drilling a lot of {\ibf Holes}, consistently, both at the (thickened) levels of $X^2$ and $fX^2$. By ``Holes'' we mean here deletions which, to be repaired back, need additions of handles of index $\lambda = 2$. Working now only with objects with Holes, we will be able to fill in the empty space left inside the ${\rm Ditch} \, (n)$ {\ibf only} for $1-\varepsilon_n \leq z \leq 1$ where, remember $\underset{n = \infty}{\rm lim} \, \varepsilon_n = 0$. This replaces the trouble-making (\ref{eq0.11}) by the following item
\begin{equation}
\label{eq0.15}
\lim_{n = \infty} \{\mbox{{\it truncated}} \ {\rm Ditch} \, (n)_y \} = (x = x_{\infty} \, , \ z = \varepsilon \, , \ t = 1) \, ,
\end{equation}
which is now on the good side, both of the Stallings barrier and of the non metrizability barrier. But, after this {\ibf partial ditch-filling process}, we have to go back to the Holes and put back the corresponding deleted material. This far from trivial issue will be discussed later.

\smallskip

But before really going on, I will make a parenthetical comment which some readers may find useful. There are many differences between the present work (of which this present paper is only a first part, but a long complete version exists too, in hand-written form) and my ill-fated, by now dead $\pi_1^{\infty} \tilde M^3 = 0$ attempt [Pr\'epublications Orsay 2000-20 and 2001-57]. Of course, a number of ideas from there found their way here too. In the dead papers I was also trying to mimick the zipping by inclusions; but there, this was done by a system of ``gutters'' added in $2^{\rm d}$ or $3^{\rm d}$, before any thickening into high dimensions. Those gutters came with fatal flaws, which opened a whole Pandora's box of troubles. These turned out to be totally unmanageable, short of some input of new ideas. And, by the time a first whiff of such ideas started popping up, Perelman's announcement was out too. So, I dropped then the whole thing, turning to more urgent tasks.

\smallskip

It is only a number of years later that this present work grew out of the shambles. By contrast with the low-dimensional gutters from the old discarded paper, the present ditches take full advantages of the {\ibf additional} dimensions; I use here the word ``additional'', by opposition to the mere high dimensions. The spectre of non metrizability which came with the ditches, asked then for Holes the compensating curves of which are dragged all over the place by the inverse of the zipping flow, far from their normal location, {\it a.s.o.} End of practice.

\smallskip

With the Holes in the picture, the $\Theta^n (X^2)$, $S_u (\tilde M^3 (\Gamma))$, will be replaced by the highly non-simply-connected smooth $n$-manifolds with non-empty boundary $\Theta^n (X^2) - H$ and $S_u (\tilde M^3 (\Gamma)) - H$. Here the ``$-H$'' stands for ``with the Holes having been drilled'' or, in a more precise language, as it will turn out, with the 2-handles which correspond to the Holes in question, deleted.

\smallskip

The manifold $S_u (\tilde M^3 (\Gamma)) - H$ comes naturally equipped with a PROPER framed link
\begin{equation}
\label{eq0.16}
\sum_1^{\infty} \, C_n \overset{\alpha}{\longrightarrow} \partial \, (S_u (\tilde M^3 (\Gamma)) - H) \, ,
\end{equation}
where ``$C$'' stands for ``curve''. This framed link is such that, when one adds the corresponding 2-handles to $S_u (\tilde M^3 (\Gamma)) - H$, then one gets back the $S_u (\tilde M^3 (\Gamma))$.

\smallskip

This was the $S_u$ story with holes, which is quite simple-minded, and we move now on to $S_b$. One uses now the $\Theta^n (X^2) - H$, {\it i.e.} the thickened $X^2$, with Holes, as a starting point for the geometric realization of the zipping process, rather than starting from the $\Theta^n (X^2)$ itself. The Holes allow us to make use of a partial ditch filling, {\it i.e.} to use now the truncation $1-\varepsilon_n \leq z \leq 1$, which was already mentioned before. This has the virtue of putting us on the good side of all the various barriers which we have to respect. The end-product of this process is another smooth $n$-dimensional manifold, which we call $S_b (\tilde M^3 (\Gamma)) - H$. This comes with a relatively easy diffeomorphism
$$
\xymatrix{
S_b (\tilde M^3 (\Gamma)) - H \ar[rr]_{\approx}^{\eta}   &&S_u (\tilde M^3 (\Gamma)) - H \, . 
} 
$$
Notice that only ``$S_b (\tilde M^3 (\Gamma)) - H$'' is defined, so far, and not yet the full $S_b (\tilde M^3 (\Gamma))$ itself. Anyway here comes the following fact, which is far from being trivial.

\bigskip

\noindent {\bf Lemma 1.2.} {\it There is a second PROPER framed link
\begin{equation}
\label{eq0.17}
\sum_1^{\infty} \, C_n \overset{\beta}{\longrightarrow} \partial \, (S_b (\tilde M^3 (\Gamma)) - H) 
\end{equation}
which has the following two features.}

\medskip

1) {\it The following diagram is commutative, {\ibf up to a homotopy}}
\begin{equation}
\label{eq0.18}
\xymatrix{
S_b (\tilde M^3 (\Gamma)) -H \ar[rr]_{\eta}   &&S_u (\tilde M^3 (\Gamma)) - H   \\ 
&\overset{\infty}{\underset{1}{\sum}} \, C_n \ar[ul]^{\beta} \ar[ur]_{\alpha}
} 
\end{equation}
{\it The homotopy above, which is {\ibf not} claimed to be PROPER, is compatible with the framings of $\alpha$ and $\beta$.}

\medskip

2) {\it We {\rm define} now the smooth $n$-dimensional manifold
$$
S_b (\tilde M^3 (\Gamma)) = (S_b (\tilde M^3 (\Gamma)) - H) + \{\mbox{\rm the $2$-handles which are defined by} 
$$
\begin{equation}
\label{eq0.19}
\mbox{\rm the framed link $\beta$ {\rm (1.17)}}\} \, .
\end{equation}
It is claimed that this manifold $S_b (\tilde M^3 (\Gamma))$ is} GSC.

\bigskip

Without even trying to prove anything, let us just discuss here some of the issues which are involved in this last lemma.

\smallskip

In order to be able to discuss the (\ref{eq0.17}), let us look at a second toy-model, the next to appear, in increasing order of difficulty, after the one already discussed, when formulae (\ref{eq0.9}) to (\ref{eq0.12}) have been considered. We keep now the same $\underset{1}{\overset{\infty}{\sum}} \, V_n$, and just replace the former $W$ by $W_1 = (y=0) \cup (z=0)$. The $M$ from (\ref{eq0.9}) becomes now the following non-compact 3-manifold with boundary
$$
M_1 = [(-\varepsilon \leq y \leq \varepsilon) \cup (-\varepsilon \leq z \leq \varepsilon) - \{\mbox{the present contribution of} 
$$
\begin{equation}
\label{eq0.20}
{\rm LIM} \, M_2 (f) \}] \cup \sum_{1}^{\infty} \, V_n \times (x_n - \varepsilon_n \leq x \leq x + \varepsilon_n ) \, ;
\end{equation}
the reader should not find it hard to make explicit the contribution of ${\rm LIM} \, M_2 (f)$ here. Also, because we have considered only $(y=0) \, \cup \, (z=0)$ and not the slightly more complicated disjoint union $(y=0) + (z=0)$, which comes with triple points, there is still no difference, so far, between ${\rm LIM} \, M_2 (f)$ and $f \, {\rm LIM} \, M_2 (f)$.

\smallskip

When we have discussed the previous local model, then the
$$
{\rm DITCH} \, (n) \underset{\rm def}{=} \ \bigcup_y \ {\rm DITCH} \, (n)_y
$$
was concentrated in the neighbourhood of the rectangle
$$
(x = x_n \, , \ -\infty < y < \infty \, , \ - \varepsilon \leq z \leq \varepsilon \, , \ t=1 ) \, .
$$
Similarly, the present ${\rm DITCH} \, (n)$ will be concentrated in a neighbourhood of the $2$-dimensional infinite cross
$$
(x = x_n \, , \ (-\varepsilon \leq y \leq \varepsilon) \cup (-\varepsilon \leq z \leq \varepsilon) \, , \ t=1) \, .
$$
It is only the $V_n$'s which see Holes. Specifically, $V_n - H$ is a very thin neighbourhood of the $1^{\rm d}$ cross $(y = + \varepsilon) \cup (z = +\varepsilon)$, living at $x=x_n$, for some fixed $t$, and fitting inside ${\rm DITCH} \, (n)$, without touching anything at the level of the $\{ W_1$ thickened in the high dimension, and with the DITCH deleted$\}$. But it does touch to four Holes, corresponding to the four corners. The action takes place in the neighbourhood of $t=1$, making again full use of the additional dimensions, supplementary to those of $M_1$ (\ref{eq0.20}).

\smallskip

With this set-up, when we try to give a ``normal'' definition for the link $\beta$ in (\ref{eq0.17}), then we encounter the following kind of difficulty, and things only become worse when triple points of $f$ are present too. Our normal definition of $\beta \, C_n$ (where $C_n$ is here the generic boundary of one of our four Holes), is bound to make use of arcs like $I_n = (x=x_n \, , \ y = +\varepsilon \, , \ -\varepsilon \leq z \leq \varepsilon \, , \ t = {\rm const})$, or $I_n = (x = x_n \, , \ - \varepsilon \leq y \leq \varepsilon \, , \ z = + \varepsilon \, , \ t = {\rm const})$ which, in whatever $S_b (\tilde M^3 (\Gamma)) - H$ may turn out to be, accumulate at finite distance. So, the ``normal'' definition of $\beta$ fails to be PROPER, and here is what we will do about this. Our arcs $I_n$ come naturally with (not completely uniquely defined) double points living in $\Psi (f) = \Phi (f)$, attached to them and these have zipping paths, like in the Representation Theorem. The idea is then to push the arcs $I_n$ back, along (the inverses of) the zipping paths, all the way to the singularities of $f$, and these do not accumulate at finite distance. In more precise terms, the arcs via which we replace the $I_n$'s inside $\beta \, C_n$ come closer and closer to $f \, {\rm LIM} \, M_2 (f)$ as $n \to \infty$, and this last object lives at infinity. This is the {\it correct} definition of $\beta$ in (\ref{eq0.17}). The point is that, now $\beta$ is PROPER, as claimed in our lemma. But then, there is also a relatively high price to pay for this. With the roundabout way to define $\beta$ which we have outlined above, the point 2) in our lemma, {\it i.e.} the basic property $S_b \in$ GSC, is no longer the easy obvious fact which it normally should be; it has actually become non-trivial. Here is where the difficulty sits. To begin with, our $X^2$ which certainly is GSC, by construction, houses two completely distinct, not everywhere well-defined flow lines, namely the zipping flow lines and the collapsing flow lines stemming from the easy id $+$ nilpotent geometric intersection matrix of $X^2$. By itself, each of these two systems of flow lines is quite simple-minded and controlled, but not so the combined system, which can exhibit, generally speaking, closed oriented loops.

\smallskip

These {\ibf bad cycles} are in the way for the GSC property of $S_b$; they introduce unwanted terms in the corresponding geometric intersection matrix. An infinite machinery is required for handling this new problem: the bad cycles have to be, carefully, pushed all the way to infinity, out of the way. Let me finish with Lemma~1.2 by adding now the following item. Punctures are normally achieved by infinite sequences of dilatations, but if we locate the $\beta \, \underset{1}{\overset{\infty}{\sum}} \, C_n$ over the regions created by them, this again may introduce unwanted terms inside the geometric intersection matrix of $S_b$, making havoc of the GSC feature. In other words, we have now additional restrictions concerning the punctures, going beyond what the Stallings barrier would normally tolerate. In more practical terms, what this means that the use of punctures is quite drastically limited, by the operations which make our $\beta$ (\ref{eq0.17}) PROPER. This is as much as I will say about the Lemma~1.2, as such.

\smallskip

But now, imagine for a minute that, in its context, we would also know that (\ref{eq0.18}) commutes up to PROPER homotopy. In that hypothetical case, in view of the high dimensions which are involved, the (\ref{eq0.18}) would commute up to isotopy too. This would prove then that the manifolds $S_u (\tilde M^3 (\Gamma))$ and $S_b (\tilde M^3 (\Gamma))$ are diffeomorphic. Together with 2) in Lemma~1.2, we would then have a proof of our desired Main Lemma. But I do not know how to make such a direct approach work, and here comes what I can offer instead.

\smallskip

The starting point of the whole $S_u / S_b$ story, was an {\ibf equivariant} representation theorem for $\tilde M^3 (\Gamma)$, namely our $2^{\rm d}$ Representation Theorem and, from there on, although we have omitted to stress this until now, everything we did was supposed to be equivariant all along: the zipping, the Holes, the (\ref{eq0.16}) $+$ (\ref{eq0.17}), {\it a.s.o.} are all equivariant things. 

\smallskip

Also, all this equivariant discussion was happening upstairs, at the level of $\tilde M^3 (\Gamma)$. But, being equivariant it can happily be pushed down to the level of $M^3 (\Gamma) = \tilde M^3 (\Gamma) / \Gamma$. Downstairs too, we have now two, still non-compact manifolds, namely $S_u (M^3 (\Gamma)) \underset{\rm def}{=} S_u (\tilde M^3 (\Gamma)) / \Gamma$ and $S_b (M^3 (\Gamma)) \underset{\rm def}{=} S_b (\tilde M^3 (\Gamma)) / \Gamma$.

\smallskip

What these last formulae mean, is also that we have
\begin{equation}
\label{eq0.21}
S_u (M^3 (\Gamma))^{\sim} = S_u (\tilde M^3 (\Gamma)) \quad \mbox{and} \quad S_b (M^3 (\Gamma))^{\sim} = S_b (\tilde M^3 (\Gamma)) \, , 
\end{equation}
let us say that $S_u$ and $S_b$ are actually functors of sorts.

\smallskip

But before really developping this new line of thought, we will have to go back to the diagram (\ref{eq0.18}) which, remember, commutes up to homotopy. Here, like in the elementary text-books, we would be very happy now to change
$$
\alpha \, C_n \sim \eta \beta \, C_n \quad \mbox{into (something like)} \quad \alpha \, C_n \cdot \eta \beta \, C_n^{-1} \sim 0 \, .
$$
This is less innocent than it may look, since in order to be of any use for us, the infinite system of closed curves
\begin{equation}
\label{eq0.22}
\Lambda_n \underset{\rm def}{=} \alpha \, C_n \cdot \eta \beta \, C_n^{-1} \subset \partial \, (S_u (\tilde M^3 (\Gamma)) - H) \, , \ n = 1,2, \ldots
\end{equation}
better be PROPER (and, of course, equivariant too). The problem here is the following and, in order not to over complicate our exposition, we look again at our simplest local model. Here, in the most difficult case at least, the curve $\alpha \, C_n$ runs along $(x = x_n \, , \ z = - \varepsilon)$ while $\eta \beta \, C_n$ runs along $(x = x_n \, , \ z = + \varepsilon)$. The most simple-minded procedure for defining (\ref{eq0.22}) would then be to start by joining them along some arc of the form
$$
\lambda_n = (x = x_n \, , \ y = {\rm const} \, , \ -\varepsilon \leq z \leq \varepsilon \, , \ t = {\rm const}) \, .
$$

But then, for the very same reasons as in our previous discussion of the mutual contradictory effects of the two barriers (Stallings and non-metrizability), this procedure is certainly not PROPER. The cure for this problem is to use, once again, the same trick as for defining a PROPER $\beta$ in (\ref{eq0.17}), namely to push the stupid arc $\lambda_n$ along the (inverse of) the zipping flow, all the way back to the singularities, keeping things all the time close to $f \, {\rm LIM} \, M_2 (f)$, {\it i.e.} close to $(x=x_{\infty} \, , \ z = \pm \varepsilon \, , \ t=1)$, in the beginning at least.

\smallskip

The next lemma sums up the net result of all these things.

\bigskip

\noindent {\bf Lemma 1.3.} {\it The correctly defined system of curves {\rm (\ref{eq0.22})} has the following features}
\begin{itemize}
\item[0)] {\it It is equivariant (which, by now, does not cost much),}
\item[1)] {\it It is PROPER,}
\item[2)] {\it For each individual $\Lambda_n$, we have a null homotopy
$$
\Lambda_n \sim 0 \quad \mbox{in} \quad S_u (\tilde M^3 (\Gamma)) - H \, .
$$
[Really this is in $\partial (S_u (\tilde M^3 (\Gamma)) - H)$, but we are very cavalier now concerning the distinction between $S_u$ and $\partial S_u$; the thickening dimension is very high, anyway.]}
\item[3)] {\it As a consequence of the bounded zipping length in the Representation Theorem, our system of curves $\Lambda_n$ has {\ibf uniformly bounded length}.}
\end{itemize}

\bigskip

Lemma~1.3 has been stated in the context of $\tilde M^3 (\Gamma)$, upstairs. But then we can push it downstairs to $M^3 (\Gamma)$ too, retaining 1), 2), 3) above. So, from now on, we consider the correctly defined system $\Lambda_n$ {\it downstairs}. Notice here that, although $M^3 (\Gamma)$ is of course compact (as a consequence of $\Gamma$ being finitely presented), the $S_u (M^3 (\Gamma)) - H$ and $S_u (M^3 (\Gamma))$ are certainly not. 

\medskip

So, the analogue of 1) from Lemma~1.3, which reads now
$$
\lim_{n = \infty} \, \Lambda_n = \infty \, , \ \mbox{inside} \ S_u (M^3 (\Gamma))-H \, , 
$$
is quite meaningful. Finally, here is the

\bigskip

\noindent {\bf Lemma 1.4.} (KEY FACT) {\it The analogue of diagram {\rm (\ref{eq0.18})} downstairs, at the level of $M^3 (\Gamma)$, commutes now up to PROPER homotopy.}

\bigskip

Before discussing the proof of this key fact, let us notice that it implies that $S_u (M^3 (\Gamma)) \underset{\rm DIFF}{=} S_b (M^3 (\Gamma))$ hence, {\it via} (\ref{eq0.21}), {\it i.e.} by ``functoriality'', we also have
$$
S_u (\tilde M^3 (\Gamma)) \underset{\rm DIFF}{=} S_b (\tilde M^3 (\Gamma)) \, , 
$$
making that $S_u (\tilde M^3 (\Gamma))$ is GSC, as desired; see here 2) in Lemma~1.2 too.

\smallskip

All the rest of the discussion is now downstairs, and we will turn back to Lemma~1.4. Here, the analogue of 2) from Lemma~1.3 is, of course, valid downstairs too, which we express as follows
$$
\mbox{For every $\Lambda_n$ there is a singular disk $D_n^2 \subset S_u (M^3 (\Gamma)) - H$,} 
$$
\begin{equation}
\label{eq0.23}
\mbox{with $\partial D^2 = \Lambda_n$.} 
\end{equation}
With this, what we clearly need now for Lemma~1.4, is something like (\ref{eq0.23}), but with the additional feature that $\underset{n=\infty}{\lim} \, D_n^2 = \infty$, inside $S_u (M^3 (\Gamma)) - H$. In a drastically oversimplified form, here is how we go about this. Assume, by contradiction, that there is a compact set $K \subset \partial (S_u (M^3 (\Gamma)) - H)$ and a subsequence of $\Lambda_1 , \Lambda_2 , \ldots$, which we denote again by exactly the same letters, such that for {\ibf any} corresponding system  of singular disks cobounding it, $D_1^2 , D_2^2 , \ldots$, we should have $K \cap D_n^2 \ne \emptyset$, for all $n$'s.

\smallskip

We will  show now that this itself, leads to a contradiction. Because $M^3 (\Gamma)$ is compact ($\Gamma$ being finitely presented), we can {\ibf compactify} $S_u (M^3 (\Gamma)) - H$ by starting with the normal embedding $S_u (M^3 (\Gamma)) - H \subset M^3 (\Gamma) \times B^N$, $N$ large; inside this compact metric space, we take then the closure of $S_u (M^3 (\Gamma)) - H$. For this compactification which we denote by $(S_u (M^3 (\Gamma)) - H)^{\wedge}$ to be nice and useful for us, we have to be quite careful about the exact locations and sizes of the Holes, but the details of this are beyond the present outline. This good compactification is now
$$
(S_u (M^3 (\Gamma)) - H)^{\wedge} = (S_u (M^3 (\Gamma)) - H) \cup E_{\infty} \, , 
$$
where $E_{\infty}$ is the compact space which one has to add at the infinity of $S_u (M^3 (\Gamma))$ $- \, H$, so as to make it compact. It turns out that $E_{\infty}$ is moderately wild, failing to the locally connected, although it has plenty of continuous arcs embedded inside it.

\smallskip

Just by metric compactness, we already have
$$
\lim_{n = \infty} {\rm dist} (\Lambda_n , E_{\infty}) = 0 
$$
and, even better, once we know that the lengths of $\Lambda_n$ are uniformly bounded, there is a subsequence $\Lambda_{j_1} , \Lambda_{j_2} , \Lambda_{j_3} , \ldots$ of $\Lambda_1 , \Lambda_2 , \Lambda_3 , \ldots$ and a continuous curve $\Lambda_{\infty} \subset E_{\infty}$, such that $\Lambda_{j_1} , \Lambda_{j_2} , \ldots$ converges uniformly to $E_{\infty}$. To be pedantically precise about it, we have
\begin{equation}
\label{eq0.24}
{\rm dist} \, (\Lambda_{j_n} , \Lambda_{\infty}) = \varepsilon_n \, , \ \mbox{where} \ \varepsilon_1 > \varepsilon_2 > \ldots > 0 \ \mbox{and} \ \lim_{n = \infty} \varepsilon_n = 0 \, .
\end{equation}
Starting from this data, and injecting also a good amount of precise knowledge concerning the geometry of $S_u (M^3 (\Gamma))$ (knowledge which we actually have to our disposal, in real life), we can construct a region $N = N (\varepsilon_1 , \varepsilon_2 , \ldots) \subset S_u (M^3 (\Gamma)) - H$, which has the following features

\medskip

A) The map $\pi_1 N \longrightarrow \pi_1 (S_u (M^3 (\Gamma)) - H)$ injects;

\medskip

B) The $E_{\infty}$ lives, also, at the infinity of $N$, to which it can be glued, and there is a retraction
$$
N \cup E_{\infty}  \overset{R}{\longrightarrow} E_{\infty} \, ;  
$$

C) There is an ambient isotopy of $S_u (M^3 (\Gamma)) - H$, which brings all the $\Lambda_{j_1} , \Lambda_{j_2} , \ldots$ inside $N$. After this isotopy, we continue to have (\ref{eq0.24}), or at least something very much like it.

\smallskip

The reader may have noticed that, for our $N$ we studiously have avoided the word ``neighbourhood'', using ``region'' instead; we will come back to this.

\smallskip

Anyway, it follows from A) above that our $\Lambda_{j_1} , \Lambda_{j_2} , \ldots \subset N$ (see here C)) bound singular disks in $N$. Using B), these disks can be brought very close to $E_{\infty}$, making them disjoined from $K$. In a nutshell, this is the contradiction which proves what we want. We will end up with some comments.

\smallskip

To begin with, our $N = N (\varepsilon_1 , \varepsilon_2 , \ldots)$ is by no means a neighbourhood of infinity, it is actually too thin for that, and its complement is certainly not pre-compact.

\smallskip

In the same vein, our argument which was very impressionistically sketched above, is certainly not capable of proving things like $\pi_1^{\infty} (S_u (M^3 (\Gamma)) - H) = 0$, which we do {\ibf not} claim, anyway.

\smallskip

Finally, it would be very pleasant if we could show that $\Lambda_{\infty}$ bounds a singular disk inside $E_{\infty}$ and deduce then our desired PROPER homotopy for the $\Lambda_{i_n}$'s from this. Unfortunately, $E_{\infty}$ is too wild a set to allow such an argument to work.

\smallskip

This ends the sketch of the proof that all $\Gamma$'s are QSF and now I will present some

\bigskip

\noindent CONJECTURAL FURTHER DEVELOPMENTS. The first possible development concerns a certain very rough classification of the set of all (finitely presented) groups. We will say that a group $\Gamma$ is {\ibf easy} if it is possible to find for it some 2-dimensional representation with {\ibf closed} $M_2 (f)$. We certainly do not mean here some equivariant representation like in the main theorem of the present paper which, most likely, will have ${\rm LIM} \, M_2 (f) \ne \emptyset$. We do not ask for anything beyond GSC and $\Psi = \Phi$. On the other hand, when there is {\ibf no} representation for $\Gamma$ with a closed $M_2 (f)$, {\it i.e.} if for any 2$^{\rm d}$ representation we have ${\rm LIM} \, M_2 (f) \ne \emptyset$, then we will say that the group $\Gamma$ is {\ibf difficult}.

\smallskip

Do not give any connotations to these notions of easy group {\it versus} difficult group, beyond the tentative technical definitions given here.

\smallskip

But let us move now to 3-dimensional representations too. All the 3-dimen\-sional representations $X \overset{f}{\longrightarrow} \tilde M^3 (\Gamma)$ are such that $X$ is a union of {\ibf ``fundamental domains''}, pieces on which $f$ injects and which have sizes which are uniformly bounded ({\it i.e.} of bounded diameters). With this, $\Gamma$ will be said now to be easy if for any compact $K \subset \tilde M^3 (\Gamma)$ there is a representation $X \overset{f}{\longrightarrow} \tilde M^3 (\Gamma)$, (possibly depending on $K$), such that only finitely many fundamental domains $\Delta \subset X$ are such that $K \cap f\Delta \ne \emptyset$. There are clearly two distinct notions here, the one just stated and then also the stronger one where a same $(X^2 , f)$ is good for all $K$'s. This last one should certainly be equivalent to the 2-dimensional definition which was given first. But at the present stage of the discussion, I prefer to leave, temporarily, a certain amount of fuzzynen concerning the $3^{\rm d}$ notion of ``difficult group''. Things will get sharper below (see that CONJECTURE~1.5 which follows). Anyway, the general idea here is that the easy groups are those which manage to avoid the {\ibf Whitehead nightmare} which is explained below, and for which we also refer to \cite{17}. It was said earlier that we have a not too difficult implication $\{ \Gamma$ is easy$\} \Longrightarrow \{ \Gamma$ is QSF$\}$. I believe this holds even with the $K$-dependent version of ``easy'', but I have not checked this fact.

\smallskip

So, the difficult groups are defined now to be those for which {\ibf any} representation $(X,f)$ exhibits the following Whitehead nightmare
$$
\mbox{For any compact $K \subset \tilde M^3 (\Gamma)$ there are INFINITELY many} 
$$
$$
\mbox{fundamental domains $\Delta \subset X$ {\it s.t.} $K \cap f\Delta \ne \emptyset$.}
$$
The Whitehead nightmare above is closely related to the kind of processes {\it via} which the Whitehead manifold ${\rm Wh}^3$ itself or, even more seriously, the Casson Handles, which have played such a proeminent role in Freedman's proof of the TOP $4$-dimensional Poincar\'e Conjecture, are constructed; this is where the name comes from, to begin with.

\smallskip

Here is the story behind these notions. Years ago, various people like Andrew Casson, myself, and others, have written papers (of which \cite{18}, \cite{19}, \cite{7},$\ldots$ are only a sample), with the following general gist. It was shown, in these papers, that if for a closed 3-manifold $M^3$, the $\pi_1 \, M^3$ has some kind of nice geometrical features, then $\pi_1^{\infty} \tilde M^3 = 0$; the list of nice geometrical features in question, includes Gromov hyperbolic (or more generally almost convex), automatic (or more generally combable), {\it a.s.o.} The papers just mentioned have certainly been superseded by Perelman's proof of the geometrization conjecture, but it is still instructive to take a look at them from the present vantage point: with hindsight, what they actually did, was to show that, under their respective geometric assumptions, the $\pi_1 \, M^3$ way easy, hence QSF, hence $\pi_1^{\infty} = 0$. Of the three implications involved here the really serious one was the first and the only $3$-dimensional one, the last.

\smallskip

As a small aside, in the context of those old papers mentioned above, both Casson and myself we have developed various group theoretical concepts, out of which Brick and Mihalik eventually abstracted the notion of QSF. For instance, I considered ``Dehn exhaustibility'' which comes with something looking, superficially like the QSF, except that now both $K$ and $X$ are smooth and, more seriously so, $f$ is an {\ibf immersion}. Contrary to the QSF itself, the Dehn exhaustibility fails to be presentation independent, but I have still found it very useful now, as an ingredient for the proof of the implication (\ref{eq0.8}). Incidentally also, when groups are concerned, Daniele Otera~\cite{10} has proved that QSF and Dehn exhaustibility are equivalent, in the same kind of weak sense in which he and Funar have proved that QSF and GSC are equivalent, namely $\Gamma \in {\rm QSF}$ iff $\Gamma$ possesses {\ibf some} presentation $P$ with $\tilde P$ Dehn exhaustible and/or GSC. See here \cite{5}, \cite{10} and \cite{4}.

\smallskip

So, concerning these same old papers as above, if one forgets about three dimensions and about $\pi_1^{\infty}$ (which I believe to be essentially a red herring in these matters), what they actually prove too, between the lines, is that any $\Gamma$ which satisfies just the geometrical conditions which those papers impose, is in fact easy.

\smallskip

What then next? For a long time I have tried unsuccessfully, to prove that any $\pi_1 \, M^3$ is easy. Today I believe that this is doable, provided one makes use of the geometrization of 3-manifolds in its full glory, {\it i.e.} if one makes full use of Perelman's work. But then, rather recently, I have started looking at these things from a different angle. What the argument for the THEOREM $\forall \, \Gamma \in {\rm QSF}$ does, essentially, is to show that even if $\Gamma$ is difficult, it still is QSF. Then, I convinced myself that the argument in question can be twisted around and then used in conjunction with the fact that we already know now that $\Gamma$ is QSF, so as to prove a much stronger result, which I only state here as a conjecture, since a lot of details are still to be fully worked out. Here it is, in an improved version with respect to the earlier, related statement, to be found in \cite{24}.

\bigskip

\noindent {\bf Conjecture 1.5.} I) {\it $2^{\rm d}$ form: For every $\Gamma$ there is a $2^{\rm d}$ representation 
$$
X^2 \overset{f}{\longrightarrow} \tilde M^3 (\Gamma)
$$
such that ${\rm LIM} \, M_2 (f) = \emptyset$, i.e. such that the double points set $M_2 (f) \subset X^2$ is {\ibf closed}}

\medskip

II) {\it $3^{\rm d}$ form: For every $\Gamma$ there is a $3^{\rm d}$ representation
$$
X^3 \longrightarrow \tilde M^3 (\Gamma) \, ,
$$
coming with a function $Z_+ \overset{\mu}{\longrightarrow} Z_+$, such that for any fundamental domain $\delta \subset \tilde M^3 (\Gamma)$, there are at most $\mu (\Vert \delta \Vert)$ fundamental domains $\Delta \subset X$ such that $f \Delta \cap \delta \ne \emptyset$. Here $\Vert \delta \Vert$ is the word-length norm coming from the group $\Gamma$.

\medskip

In other words, {\ibf all groups are easy}, i.e. they all can avoid the Whitehead nightmare.}

\bigskip

This certainly implies, automatically, that $\Gamma \in {\rm QSF}$ and also, according to what I have already said above, it should not really be a new thing for $\Gamma = \pi_1 \, M^3$. But in order to get to that, one has to dig deeper into the details of the Thurston Geometrization, than one needs to do for $\pi_1 \, M^3 \in$ QSF, or at least so I think. 

\smallskip

All this potential development concerning easy and difficult groups is work in progress. The next possible developments which I will briefly review now are not even conjectures but rather hypothetical wild dreams, bordering to science-fiction.

\smallskip

To begin with, when it comes to zipping then there is the notion of COHERENCE which I will not redefine here, a very clear exposition of it can be found in \cite{6}. Coherence is largely irrelevant in dimensions $n > 4$ but then, it is tremendously significant when the dimension is four.

\bigskip

\noindent {\bf Question 1.6.} Given any $\Gamma$, can one always find some $2^{\rm d}$ representation $X^2 \overset{f}{\longrightarrow} \tilde M^3 (\Gamma)$ for which there is a COHERENT zipping strategy?

\bigskip

There should be no question here of anything like a nice, equivariant representation. My wild tentative conjecture is that one might be able to prove the statement with question mark above by using the techniques which are very briefly alluded to in \cite{22}, \cite{23}. I am in particular thinking here about the proof of the so-called COHERENCE THEOREM (also stated in \cite{6}). The proofs in question (and also for the rest of \cite{22}, \cite{23}) are actually already completely (hand-)written. If the answer to question~1.6 would be yes, then we could construct a {\ibf 4-dimensional} $S_u \tilde M^3 (\Gamma)$, most likely not GSC and certainly not equivariant. But this would come with a boundary, an open non simply connected wild $3$-manifold $V^3$, with $\pi_1^{\infty} V^3 \ne 0$. The point is that this $V^3$ without actually really supporting an action of $\Gamma$ would still be in a certain sense related to $\Gamma$. 

\smallskip

This could  conceivably have some interest.

\smallskip

Next, here are two related questions for which we would like to have an answer: How can be reconcile the statement that all $\Gamma$'s are QSF, with the commonly accepted idea that any property valid for all $\Gamma$'s has to be trivial? And then, if we belive Conjecture~1.5 above, where do the ``difficult objects'' in group theory, if any, hide? I believe that there should be a class of objects which, very tentatively I will call here ``quasi groups'' and among which the finitely presented groups should presumably live somehow like the rational numbers among the reals. My ``quasi'' refers here to quasiperiodic as opposed to periodic and, of course also, to quasi-crystals and/or to Penrose tilings. I cannot even offer here a conjectural definition for ``quasi groups''. I might have some guesses about what a quasi group presentation should be and then I can already see a serious difficulty in deciding when two presentations define the same object.

\bigskip

\centerline{* \ * \ * \ * \ *}

\medskip

The many criticisms, comments and suggestions which David Gabai has made in connection with my earlier ill-fated attempt of proving $\pi_1^{\infty} \, \tilde M^3 = 0$, were essential for the present work, which could not have existed without them. Like in other occasions too, his help was crucial for me.

\smallskip

Thanks are also due to Louis Funar and Daniele Otera for very useful conversations. Actually, it was Daniele who, at that time my PhD student, first told me about QSF, and who also insisted that I should look into it.

\smallskip

In 2007 and 2008 I have lectured on these matters in the Orsay Geometric Group Theory Seminar. Thanks are due to Frederic Harglund and to the other members of the Seminar, for helpful comments and conversations.

\smallskip

Finally, I wish to thank IH\'ES for its friendly help, C\'ecile Cheikhchoukh for the typing and Marie-Claude Vergne for the drawings.

\section{Zipping}\label{sec1}
\setcounter{equation}{0}

The presentations for finitely generated groups $\Gamma$, which we will use, will be singular compact 3-manifolds with boundary $M^3 (\Gamma)$, such that $\pi_1 M^3 (\Gamma) = \Gamma$. Here is how such a $M^3 (\Gamma)$ will be gotten. To a smooth handlebody $H$ ($=$ finite union of handles of index $\lambda \leq 1$), we will attach 2-handles via a generic immersion having as core a link projection
\begin{equation}
\label{eq1.1}
\sum_i (S_i^1 \times I_i) \overset{\alpha}{\longrightarrow} \partial H \, .
\end{equation}
We may assume each individual component to be embedded. The double points of $\alpha$ are little squares $S \subset \partial H$. These are the singularities ${\rm Sing} \, M^3 (\Gamma) \subset M^3 (\Gamma)$ and, to distinguish them from other singularities to occur later, they are called {\ibf immortal}. Without loss of generality, the immortal singularities $S$ live on the free part of $\partial$ (0-handles). There are no 3-handles for $M^3 (\Gamma)$.

\smallskip

Following rather closely \cite{14}, we will discuss now the equivalence relation forced by the singularities of a non-degenerate simplicial map $X \overset{f}{\longrightarrow} M$. In this story, $M$ could be any simplicial complex of dimension $n$ (although in real life it will just be $M = \tilde M^3 (\Gamma)$ or $M = M^3 (\Gamma)$ itself) and $X$ a countable Gromov multicomplex (where the intersection of two simplices is not just a common face, but a subcomplex). Since $f$ is non-degenerate, we will have $\dim X \leq \dim M$. We will need $X$'s which in general, will be not locally finite, and we will endow them with the {\ibf weak topology}. By definition, a singularity
\begin{equation}
\label{eq1.2}
x \in {\rm Sing} \, (f) \subset X
\end{equation}
is a point of $X$ in the neighbourhood of which $f$ fails to inject, {\it  i.e.} fails to be immersive. Do not mix the permanent, immortal singularities of $M^3 (\Gamma)$, {\it i.e.} $S \subset {\rm Sing} \, M^3 (\Gamma)$ with the unpermanent, {\ibf mortal} singularities (\ref{eq1.2}), of the map $f$ (as defined already in \cite{14}, \cite{15}, \cite{16} and \cite{6}).

\smallskip

Quite trivially, the map $f$ defines an equivalence relation $\Phi (f) \subset X \times X$, where $(x_1 , x_2) \in \Phi (f)$ means just that $fx_1 = fx_2$. An equivalence relation $R \subset \Phi (f)$ will be called {\ibf $f$-admissible}, if it fullfils the following condition

\medskip

\noindent (2.2.1) \quad Assume $\sigma_1 , \sigma_2 \subset X$ are two simplices of the same dimension, with $f \sigma_1 = f \sigma_2$. If we can find pairs of points $x \in {\rm int} \, \sigma_1$, $y \in {\rm int} \, \sigma_2$ with $fx = fy$, then the following IMPLICATION holds
$$
(x,y) \in R \Longrightarrow R \ \mbox{identifies} \ \sigma_1 \ \mbox{to} \ \sigma_2 \, .
$$

The equivalence relation $\Phi (f)$ and the double point set $M^2 (f) \subset X \times X$ are, of course, closely related, since $M^2 (f) = \Phi (f) - {\rm Diag} \, X$ and we will also introduce the following intermediary subset
\begin{equation}
\label{eq1.3}
\hat M^2 (f) = M^2 (f) \cup {\rm Diag} \, ({\rm Sing} \, (f)) \subsetneqq \Phi (f) \, .
\end{equation}
This has a natural structure of simplicial complex of which ${\rm Diag} \, ({\rm Sing} \, (f))$ is a subcomplex. We will endow it, for the time being, with the weak topology.

\smallskip

We will be interested in subsets $\tilde R \subset \hat M^2 (f)$ of the following form. Start with an arbitrary equivalence relation $R \subset \Phi (f)$, and then go to $\tilde R = R \cap \hat M^2 (f)$, when $R$ is an $f$-admissible equivalence relation; then we will say that $\tilde R$ itself is an {\ibf admissible set}. In other words the admissible sets are, exactly, the traces on $\hat M^2 (f)$ of $f$-admissible equivalence relations.
\begin{equation}
\label{eq1.4}
\mbox{A subset of $\hat M^2 (f)$ is admissible ({\it i.e.} it is an $\tilde R$) iff}
\end{equation}
$$
\mbox{it is both open and closed in the weak topology of $\hat M^2 (f)$.}
$$
Automatically, admissible sets are subcomplexes of $\hat M^2 (f)$.

\smallskip

When $R$ is an $f$-admissible equivalence relation, then $X/R$ is again a (multi) complex and its induced quotient space topology is again the weak topology. Also, we have a natural simplicial diagram
\begin{equation}
\label{eq1.5}
\xymatrix{
X \ar[rr]^f \ar[dr]^{\pi(R)} &&M \\
&X/R \ar[ur]_{f_1(R)}
}
\end{equation}
The $\tilde R \cap {\rm Diag} \, ({\rm Sing} \, (f)) \subset \tilde R \subset \hat M^2 (f)$ is a subcomplex, naturally isomorphic to a piece, denoted $\tilde R \cap {\rm Sing} \, (f) \subset {\rm Sing} \, (f) \subset X$, which is both open and closed in ${\rm Sing} \, (f)$.

\bigskip

\noindent {\bf Claim 2.5.1.} In the context of (\ref{eq1.5}) we have the following equality
$$
{\rm Sing} \, (f_1 (R)) = \pi (R) \, ({\rm Sing} \, (f) - \tilde R \cap {\rm Sing} \, (f)) \, .
$$
All these things, just like the next two lemmas are easy extensions to our present singular set-up, of the little theory developed in \cite{14} in the context of smooth 3-manifolds (immortal singularities were absent in \cite{14}). In fact this little, abstract non sense type theory, developed in \cite{14}, is quite general, making such extensions painless. We certainly do not even aim at full generality here.

\bigskip

\noindent {\bf Lemma 2.1.} {\it There is a {\ibf unique} $f$-admissible equivalence relation $\Psi (f) \subset \Phi (f)$ which has the following properties, which also characterize it.}

\smallskip

I) {\it When one goes to the natural commutative diagram, on the lines of} (\ref{eq1.5}), {\it i.e. to
$$
\xymatrix{
X \ar[rr]^f \ar[dr]_{\pi = \pi (\Psi (f))} &&M \\
&X/\Psi (f) \ar[ur]_{f_1= f_1 (\Psi (f))}
} 
\eqno (2.5.2)
$$
then we have ${\rm Sing} \, (f_1) = \emptyset$, i.e. $f_1$ is an immersion.}

\smallskip

II) {\it Let now $R$ be {\ibf any} equivalence relation such that $R \subset \Phi (f)$. For such an $R$ there is always a diagram like} (\ref{eq1.5}), {\it except that it may no longer be simplicial. But let us assume now also, that ${\rm Sing} \, (f_1 (R)) = \emptyset$, and also that $R \subset \Psi (f)$. Then we necessarily have $R = \Psi (f)$.}

\bigskip

We may rephrase the II) above, by saying that $\Psi (f)$ is the smallest equivalence relation, compatible with $f$, which kills all the (mortal) singularities.

\smallskip

Like in \cite{14}, we start by outlining a very formal definition for $\Psi (f)$, with which the statement above can be proved easily. For that purpose, we will introduce a new, non-Hausdorff topology for $\hat M^2 (f)$, which we will call the {\ibf $Z$-topology}. The closed sets for the $Z$-topology are the finite unions of admissible subsets. For any subset $E \subset \hat M^2 (f)$, we will denote by $C\ell_Z (E)$, respectively by $\widehat{C\ell}_Z (E) \supset C\ell_Z (E)$, the closure of $E$ in the $Z$-topology, respectively the smallest subset of $\hat M^2 (f)$ which contains $E$ and which is both $Z$-closed and admissible, {\it i.e.} the smallest equivalence relation containing $E$ and which is also $Z$-closed. One can check that the {\ibf irreducible} closed subsets are now exactly the $C\ell_Z (x,y)$ with $(x,y) \in \hat M^2 (f)$ and, for them we also have that $\widehat{C\ell}_Z (x,y) = C\ell_Z (x,y)$. With all these things, I will take for $\Psi (f)$ the following definition
$$
\Psi (f) \underset{\rm def}{=} \widehat{C\ell}_Z ({\rm Diag} \, ({\rm Sing} \, (f))) \cup {\rm Diag} \, (X) \subset \Phi (f) \, . \eqno (2.5.3)
$$
Here, the $\widehat{C\ell}_Z$ might a priori miss some of ${\rm Diag} \, (X) - {\rm Diag} \, ({\rm Sing} \, (f))$, reason for adding the full ${\rm Diag} \, (X)$ into our definition. The formula (2.5.3) defines, automatically, an $f$-admissible equivalence relation, for which we have
$$
\widetilde{\Psi (f)} = \widehat{C\ell}_Z ({\rm Diag} \, ({\rm Sing} \, (f))) \, .
$$
Also, starting from (2.5.3), it is easy to prove lemma~2.1, proceeding on the following lines, just like in \cite{14}.

\smallskip

With the notations already used in the context of the (2.5.1), one checks first that $\widetilde{\Psi (f)} \cap {\rm Sing} \, (f) = {\rm Sing} \, (f)$ and, by the same claim (2.5.1) this implies that ${\rm Sing} \, (f_1 (\Psi (f)) = \emptyset$, proving thereby I).

\smallskip

So, consider now the a priori quite arbitrary $R$ from II). We claim that $R$ has to be $f$-admissible. The argument goes as follows. The map $f_1 (R)$ being immersive, $X/R$ is Hausdorff, hence $R \subset X \times X$ is closed and hence so is also $\tilde R = X \cap \hat M^2 (f)$. Weak topology is meant here all along, in this little argument, and not the $Z$-topology.

\smallskip

Next, any non-interior point for $\tilde R$ would be a singularity for $f_1 (R)$. This implies that $\tilde R \subset \hat M^2 (f)$ is open. The offshot is that $\tilde R$ is both open and closed in the weak topology. Hence, by (\ref{eq1.4}) our $R$ has to be $f$-admissible, as it was claimed.

\smallskip

So, assume now that $R$ is $f$-admissible, with ${\rm Sing} \, f_1 (R) = \emptyset$ and with $R \subset \Psi$. By (2.5.1), $\tilde R \supset {\rm Diag} \, ({\rm Sing} \, (f))$. Being admissible, $\tilde R$ is $Z$-closed and, by assumption, it is an equivalence relation too. Being already contained in $\Psi$ it cannot fail to be equal to it. This proves our lemma~2.1. $\Box$

\bigskip

But then, once we know that our unique $\Psi (f)$ is well-defined and exists, we can also proceed now differently. I will describe now a process, called ZIPPING, which is a more direct constructive approach to $\Psi (f)$, to be used plenty in this paper.

\smallskip

So, let us look for a most efficient minimal way to kill all the mortal singularities. Start with two simplexes $\sigma_1 , \sigma_2 \subset X$ with same dimensions and $f \sigma_1 = f \sigma_2$, for which there is a singularity
$$
{\rm Sing} \, (f) \ni x_1 \in \sigma_1 \cap \sigma_2 \, .
$$
We go then to a first quotient of $X$, call it $X_1 \overset{f_1}{\longrightarrow} M$, which kills $x_1$ via a {\ibf folding map}, identifying $\sigma_1$ to $\sigma_2$. Next, start with some $x_2 \in {\rm Sing} \, (f_2)$, then repeat the same process, a.s.o. Provided things do not stop at any finite time, we get an increasing sequence of equivalence relations
$$
\rho_1 \subset \rho_2 \subset \ldots \subset \rho_n \subset \rho_{n+1} \subset \ldots \subset \Phi (f) \, . \eqno (2.5.4)
$$
The union $\rho_{\omega} = \overset{\infty}{\underset{1}{\bigcup}} \ \rho_i$ is again an equivalence relation, subcomplex of $\Phi (f)$, {\it i.e.} closed in the weak topology. The map
$$
X_{\omega} = X / \rho_{\omega} \overset{f_{\omega}}{\longrightarrow} M^3 (\Gamma)
$$
is simplicial. None of the $\rho_1 , \rho_2 , \ldots$ is $f$-admissible, of course, and, in general, neither is $\rho_{\omega}$. We then pick up some $x_{\omega} \in {\rm Sing} \, (f_{\omega})$ and go to $X/\rho_{\omega + 1} \overset{f_{\omega + 1}}{-\!\!-\!\!\!\longrightarrow} M$. From here on, one proceeds by transfinite induction and, since $X$ is assumed to be countable, the process has to stop at some countable ordinal $\omega_1$. The following things happen at this point.

\medskip

\noindent (2.5.5) \quad Using lemma~2.1 one can show that $\rho_{\omega_1} = \Psi (f)$. This makes $\rho_{\omega_1}$ canonical, but not $\omega_1$ itself. There is no unique strategy leading to $\Psi (f)$, this way, or any other way.

\bigskip

\noindent {\bf Claim 2.5.6.} One can {\it chose} the sequence (2.5.4) such taht $\omega_1 =  \omega$, {\it i.e.} such that, just with (2.5.4) we get $\rho_{\omega} = \Psi (f)$. It is this kind of sequence (2.5.4) (which is by no means unique either), which will be by definition, a {\ibf strategy for zipping $f$}, or just a ``zipping''.

\bigskip

\noindent {\bf Lemma 2.2.} {\it The following, induced map, is surjective}
$$
\pi_1 X \longrightarrow \pi_1 (X/\Psi (f)) \, .
$$

\medskip

Given the $M^3 (\Gamma)$, we will apply, in this series of papers, the little $\Psi / \Phi$ theory above in a sequence of contexts of increasing complexity. The prototype is to be shown next, and the reader should compare it with \cite{16}, \cite{18}, \cite{19}.

\smallskip

Our $M^3 (\Gamma)$, to which we come back now, is naturally divided into handles $h^{\lambda}$ of index $\lambda \leq 2$. For each such $h$ we will distinguish some individually embedded 2-cells $\varphi \subset \partial h$ called {\ibf faces}. The idea is that each face is shared exactly by two $h$'s, the whole $M^3 (\Gamma)$ being the union of the $h$'s along common faces, with the immortal singularities being created automatically, for free so to say, in the process. To begin with the attaching zone of an $h^1$ consists of two faces, little discs occuring again on the $h^0$'s. This first batch of faces, creates a system of disjoined curves $\gamma \subset \partial H$ where $H$ is the handlebody put together from the $h^0$'s and $h^1$'s; the $\gamma$'s are the boundaries of the small discs of type $h^0 \cap h^1$. Each $h^2$ has an attaching zone which, without loss of generality we may assume embedded in $\partial H$; the $\gamma$'s cut this attaching zone into long rectangles. The rectangles are a second batch of faces, the first one were discs. Each rectangle is contained either in the lateral surface of an $h^1$, going parallel to the core of $h^1$, or in the boundary of an $h^0$, connecting two discs. Two rectangles occuring on the same $\partial h^0$ may cut transversally through each other, along a little square, which is now an immortal singularity $S$. Without any loss of generality it may be assumed that all our immortal singularities $S$ occur on the lateral surface of the $0$-handles. Retain that the $S$'s are NOT faces, only parts of such.

\smallskip

Forgetting now about the Morse index $\lambda$, consider all the handles of $M^3 (\Gamma)$, call them $h_1 , h_2 , \ldots , h_p$. Some ``initial''
$$
h \in \{ h_1 , h_2 , \ldots , h_p \}
$$
will be fixed once and for all. If ${\mathcal F} (h_i)$ is the set of all the faces occuring on the boundary of $h_i$ (discs or rectangles), we note that the set ${\mathcal F} \underset{\rm def}{=} \underset{i}{\sum}  \ {\mathcal F} (h_i)$ which is of even cardinality, comes equipped with a fixed-point free involution ${\mathcal F} \overset{j}{\longrightarrow} {\mathcal F}$, with the feature that, if $x \in {\mathcal F} (h_{\ell})$ then $jx \in {\mathcal F} (h_k)$ with $k \ne \ell$, and such that
\begin{equation}
\label{eq1.6}
\mbox{Whenever $x \in {\mathcal F} (h_{\ell})$ and $jx \in {\mathcal F} (h_k)$, one glues $h_{\ell}$ to $h_k$}
\end{equation}
$$
\mbox{along $x = jx$, so as to get $M^3 (\Gamma)$.}
$$
With these things, starting from our arbitrarily chosen $h$, we will introduce a class of 3-dimensional thick paths, formally modelled on $[0,\infty)$, which are randomly exploring through $M^3 (\Gamma)$, constructed according to the following recipee.

\medskip

\noindent (2.6.1) \quad For the initial $h$, chose some $x_1 \in {\mathcal F} (h)$. There is exactly one handle in $\{ h_1 , h_2 , \ldots , h_p \} - \{ h \}$, call it $x_1 h$, such that $x_1 h$ houses $jx_1$.

\medskip

\noindent (2.6.2) \quad Inside ${\mathcal F} (x_1 h) - \{ jx_1 \}$ chose some face $x_2 \in {\mathcal F} (x_1 h)$. There is then exactly on handle in $\{ h_1 , h_2 , \ldots , h_p \} - \{ x_1 h \}$, call it $x_1 \, x_2 \, h$, such that $x_1 \, x_2 \, h$ houses $jx_2$.

$$\ldots \ldots \ldots \ldots \ldots \ldots \ldots \ldots \ldots \ldots \ldots \ldots$$

\bigskip

This kind of process continues indefinitely, producing infinitely many sequences of words written with the letters $\{ h ; x_1 , x_2 , \ldots \}$, and taking the form
$$
h , \ x_1 h , \ x_1 \, x_2 \, h , \ x_1 \, x_2 \, x_3 \, h , \ \ldots \eqno (S_{\infty})
$$

Of course, $(S_{\infty})$ can also be thought of as just an infinite word $x_1 \, x_2 \, x_3 \ldots$, written with letters $x_{\ell} \in {\mathcal F}$. The rule for constructing it makes it automatically reduced, i.e. $x_{\ell + 1} \ne jx_{\ell}$, $\forall \, \ell$.

\smallskip

{\ibf All} the sequences $(S_{\infty})$ can be put together into the following infinite, non locally finite complex, endowed with a tautological map $F$ into $M^3 (\Gamma)$, namely
\begin{equation}
\label{eq1.7}
X \underset{\rm def}{=} h \ \cup \sum_{x_1 \in {\mathcal F} (h)} x_1 h \ \cup \sum_{x_2 \in {\mathcal F} (x_1 h)} x_1 \, x_2 \, h \cup \ldots \overset{F}{\longrightarrow} M^3 (\Gamma) \, .
\end{equation}

Before we apply our little $\Phi / \Psi$ theory to this (\ref{eq1.7}), let us make a few remarks concerning it

\medskip

\noindent (2.7.1) \quad Our tree-like $X$ is certainly {\ibf arborescent}, {\it i.e.} gettable from a point by a sequence of Whitehead dilatations; and arborescence implies GSC.

\medskip

\noindent (2.7.2) \quad Up to a point the construction of $X$ is modelled on the Cayley graph. BUT there is no group involved in our construction and hence no group action either. Not all the infinite words $x_1 \, x_2 \ldots$ are acceptable for $(S_{\infty})$. There is not even a monoid present. All this makes the present construction be {\ibf not} {\it quite} like in \cite{16} or \cite{18} or \cite{19}.

\bigskip

\noindent {\bf Lemma 2.3.} {\it In the style of} (2.5.2), {\it we consider now the diagram}
\begin{equation}
\label{eq1.8}
\xymatrix{
X \ar[rr]_F \ar[dr] &&M^3 (\Gamma) \\
&X/\Psi (F) \ar[ur]_{F_1}
}
\end{equation}
{\it Then, the $X/\Psi (F) \overset{F_1}{\longrightarrow} M^3 (\Gamma)$ IS the universal covering space $\tilde M^3 (\Gamma) \overset{\pi}{\longrightarrow} M^3 (\Gamma)$.}

\bigskip

\noindent {\bf Proof.} Lemma~2.2 tells us that $\pi_1 (X/\Psi (F_1)) = 0$. Then, $F_1$ clearly has the path lifting property and it is also \'etale. End of the argument.

\smallskip

I actually like to think of this lemma~2.3 as being a sort of ``naive theory of the universal covering space''.

\smallskip

We move now to the following natural tesselation
\begin{equation}
\label{eq1.9}
\tilde M^3 (\Gamma) = \bigcup_{y \in \Gamma} \ \sum_{1}^{p} \, y \, h_i \, .
\end{equation}
Any lift of the initial $h \subset X$ to $\tilde M^3 (\Gamma)$, automatically comes with a canonical lift of the whole of $X$, which respects the handle-identities. We will denote by $f$ this lift of $F$ to $\tilde M^3 (\Gamma)$
$$
\xymatrix{
X \ar[rr]^f \ar[dr]_F &&\tilde M^3 (\Gamma) \ar[dl]^{\pi} \\
&M^3 (\Gamma)
} \eqno (2.9.1)
$$

The existence of $f$ follows from the path lifting properties of $\pi$, combined with the fact that, locally our $X$ is always something of the following form, just like in $M^3 (\Gamma)$ and/or in $\tilde M^3 (\Gamma)$
$$
\underset{\overbrace{{\mathcal F} (x_1 \ldots x_i h) \ni x_{i+1} = jx_{i+1} \in {\mathcal F} (x_1 \ldots x_{i+1} h)}}{\quad \quad \ \ x_1 \, x_2 \ldots x_i \, h \cup x_1 \, x_2 \ldots x_i \, x_{i+1} \, h \, .} 
$$
Clearly we have ${\rm Sing} \, (f) = {\rm Sing} \, (F) \subset X$ and these mortal singularities are exactly the disjoined arcs $\sigma$ of the form $\sigma = h^0 \cap h^1 \cap h^2$ where, at the level of $M^3 (\Gamma)$ the $h^{\lambda}$ and $h^{\mu}$ have in common the faces $\varphi_{\lambda \mu}$, {\it i.e.}
$$
\varphi_{01} = h^0 \cap h^1 \, , \ \varphi_{12} = h^1 \cap h^2 \, , \ \varphi_{20} = h^2 \cap h^0 \, , \ \mbox{with} \ \sigma = \varphi_{01} \cap \varphi_{12} \cap \varphi_{20} \, . \eqno (2.9.2)
$$

The $\sigma$'s, which are all mortal, live far from the immortal $S$'s. Around $\sigma$, the (\ref{eq1.7}) looks, locally, like
$$
\{\mbox{The Riemann surface of} \ \log z \} \times R \, .
$$
When one goes from $X$ to $fX$ or to $FX$, then the following things happen, as far as the immortal singularities are concerned:

\medskip

I) Two distinct immortal singularities, together with their neighbourhoods inside $X$, call them $(V_1 , S_1) \subset (X , {\rm Sing} \, X) \supset (V_2 , S_2)$, may get identified, $(V_1 , S_1) = (V_2 , S_2)$.

\medskip

II) The zipping flow of $\Psi (f) = \Psi (F)$ (see here the next lemma~2.4 too) may create new immortal singularities of $fX = \tilde M^3 (\Gamma)$, by forcing glueings of the type $h^0 \underset{\varphi}{\cup} h_i^2 \underset{\varphi}{\cup} h_j^2$. Their images via $\pi$ are then immortal singularities of $M^3 (\Gamma)$. Our $3^{\rm d}$ representation space $X$ can, itself, have immortal singularities.

\smallskip

From lemma~2.3 one can easily deduce the following

\bigskip

\noindent {\bf Lemma 2.4.} 1) {\it We have $\Psi (f) = \Psi (F) \subset X \times X$.}

\smallskip

\noindent 2) {\it The map $X \overset{f}{\longrightarrow} \tilde M^3 (\Gamma)$ (see} (2.9.1){\it ) is such that}
\begin{equation}
\label{eq1.10}
\Psi (f) = \Phi (f) \, .
\end{equation}

\medskip

Since clearly also, $X$ being collapsible is certainly GSC and the map $f$ is surjective, our $X \overset{f}{\longrightarrow} \tilde M^3 (\Gamma)$ IS a representation for $\Gamma$, albeit one which has none of the desirable features of the REPRESENTATION THEOREM. But this, presumably simplest possible representation for $\Gamma$, will be the first step towards the theorem in question.

\smallskip

But before anything else, a small drawback of our newly found representation will have to be corrected. The problem here is that, although it is a union of handles, our $X$ is {\ibf not} a {\ibf ``handlebody''}. We distinguish here HANDLEBODIES (singular of course) as being unions of handles with nice attaching maps of the $\lambda$-handles to the $(\lambda - 1)$-skeleton. We will not elaborate this notion in the most general case but, in our specific situation with only handles of indices $\lambda \leq 2$, what we demand is that, besides each 1-handle being normally attached to two 0-handles, each attaching zone of a 2-handle should find a necklace of successive handles of index $\lambda = 0$ and $\lambda = 1$, into the surface of which it should embed. This necklace itself could have repetitions, of course. Our $X$ (\ref{eq1.7}) does {\ibf not} fulfill this last condition. We will get around this difficulty by a REDEFINITION of $(X , F , M^3 (\Gamma))$, and hence of $(f,\tilde M^3 (\Gamma))$ too. More manageable mortal singularities will be gotten in this process too.

\smallskip

We proceed as follows. Start by forgetting the handlebody structure of $M^3 (\Gamma)$ and replace $M^3 (\Gamma)$ by a $2^{\rm d}$ cell-complex, via the following two stages. We first go from the handlebody $H$ (see the beginning of the section) to
$$
A^2 \underset{\rm def}{=} \{\mbox{the union of the {\it boundaries} of the various handles of index $\lambda = 0$ and $1$}\}.
$$
This is  taken to be a very finely subdivided simplicial complex, with the corresponding faces $\varphi$ occuring as subcomplexes. Our $A^2$ comes equipped with individually embedded simplicial closed curves $c_1 , c_2 , \ldots , c_{\mu}$ which correspond to the cores of the attaching zones of the $2$-handles. The immortal singularities correspond to intersections $c_i \cap c_j$ ($i \ne j$), at vertices of $A^2$. Next, we go to the cell-complex $B^2$ gotten by attaching a 2-cell $D_i^2$ to $A^2$ along each $c_i$. The $B^2$ continues to be a presentation for $\Gamma$, of course. To each of the building blocs $h , x_i , x_1 \ldots x_i \, h$ corresponds a subcomplex of $B^2$. We can glue them together by exactly the same recipee as in (\ref{eq1.7}) and, this way, we generate a purely 2-dimensional version of (\ref{eq1.7}), in the realm of $2^{\rm d}$ cell-complexes (or even simplicial complexes if we bother to subdivide a bit, see here also the lemma~2.5 below)
\begin{equation}
\label{eq1.11}
X ({\rm provisional}) \overset{F ({\rm provisional})}{-\!\!\!-\!\!\!-\!\!\!-\!\!\!-\!\!\!-\!\!\!-\!\!\!-\!\!\!-\!\!\!\longrightarrow} B^2 \, .
\end{equation}
The next stage is to get 3-dimensional again, by applying to the (\ref{eq1.11}) the standard recipee
$$
\{\mbox{cell (or simplex) of dimension} \ \lambda\} \Longrightarrow \{3^{\rm d} \ \mbox{handle of index} \ \lambda\} \, . \eqno (2.11.1)
$$
Notice that, when a cell-{\ibf complex} is changed into a union of handles via the recipee (2.11.1), then this is, automatically, a {\ibf handlebody}, and not a new union of handles. This changes (\ref{eq1.11}) into a new version of (\ref{eq1.7}), sharing all the good features (2.7.1) to (\ref{eq1.10}), which the (\ref{eq1.7}) already had, except that (2.9.2) and the $\log z$-type structure of the mortal singularities will no longer be with us (see below). This new version
\begin{equation}
\label{eq1.12}
{\rm new} \, X \overset{{\rm new} \, F}{-\!\!\!-\!\!\!-\!\!\!-\!\!\!-\!\!\!\longrightarrow} {\rm new} \, M^3 (\Gamma)
\end{equation}
will replace from now on (\ref{eq1.7}). We will have an handlebody decomposition
\begin{equation}
\label{eq1.13}
{\rm new} \, \tilde M^3 (\Gamma) = ({\rm new} \, M^3 (\Gamma))^{\sim} = \bigcup_{0 \leq \lambda \leq 2} h_i^{\lambda} \, .
\end{equation}
By construction, the new $X$ is now a singular handlebody
\begin{equation}
\label{eq1.14}
{\rm new} \, X = \bigcup_{\overbrace{\lambda , i , \alpha}} h_i^{\lambda} (\alpha) \quad \mbox{where} \quad ({\rm new} \, f)(h_i^{\lambda} (\alpha)) = h_i^{\lambda} \, .
\end{equation}

Our $h_i^{\lambda} (\alpha)$'s are here bona fide $3^{\rm d}$ handles of index $\lambda$. The $\{ \alpha \}$ is a countable system of indices, which is $(\lambda , i)$-dependent.

\bigskip

\noindent AN IMPORTANT CHANGE OF NOTATION. From now on $M^3 (\Gamma)$, $\tilde M^3 (\Gamma)$, $X$, will mean the {\ibf new} objects which we have just introduced. When the others from before, may still need to be mentioned, they will be referred to as the being {\ibf old} ones, when necessary.

\smallskip

Without any loss of generality, the immortal singularities $S$ are again corralled on the free part of the lateral surface of the $0$-handles and, also, all the desirable features from the old context continue to be with us. The new context will also have, among others, the virtue that we will be able to plug into it, with hardly any change, the theory developed in the previous paper \cite{26}.

\bigskip

\noindent THE SINGULARITIES OF THE $({\rm new}) \, X$ AND $f$. Like in the context of the singular handlebody $T \overset{F}{\longrightarrow} \tilde M^3$ from (2.9) in \cite{26}, the singularities of $f$ are now no longer modelled on $\log z$, like it was the case for the old $f$, but they occur now as follows.

\smallskip

(Description of ${\rm Sing} \, (f) = {\rm Sing} \, (F) \subset X$.) The ${\rm Sing} \, (f)$ ($=$ mortal singularities of $f$) is a union of 2-cells
\begin{equation}
\label{eq1.15}
D^2 \subset \delta \, h_{i_1}^{\lambda} (\alpha_1) \cap \partial \, h_{i_2}^{\mu} (\alpha_2) \, , \quad \mu > \lambda
\end{equation}
where for any given $\{ D^2 , (i_1 , \alpha_1 )\}$ we have {\ibf infinitely} many distinct $(i_2 , \alpha_2)$'s. Also, for any handle $h$ our notation is $\partial =$ attaching zone and $\delta =$ lateral surface. One virtue of these mortal singularities, among others is that they are amenable to the treatment from \cite{26}; the older $\log z$ type singularities were not. It should be stressed that (\ref{eq1.15}) lives far from ${\rm Sing} \, (X)$.

\smallskip

In a toy-model version, figure 2.1 offers a schematical view of the passage from the odd context to the new one.
$$
\includegraphics[width=11cm]{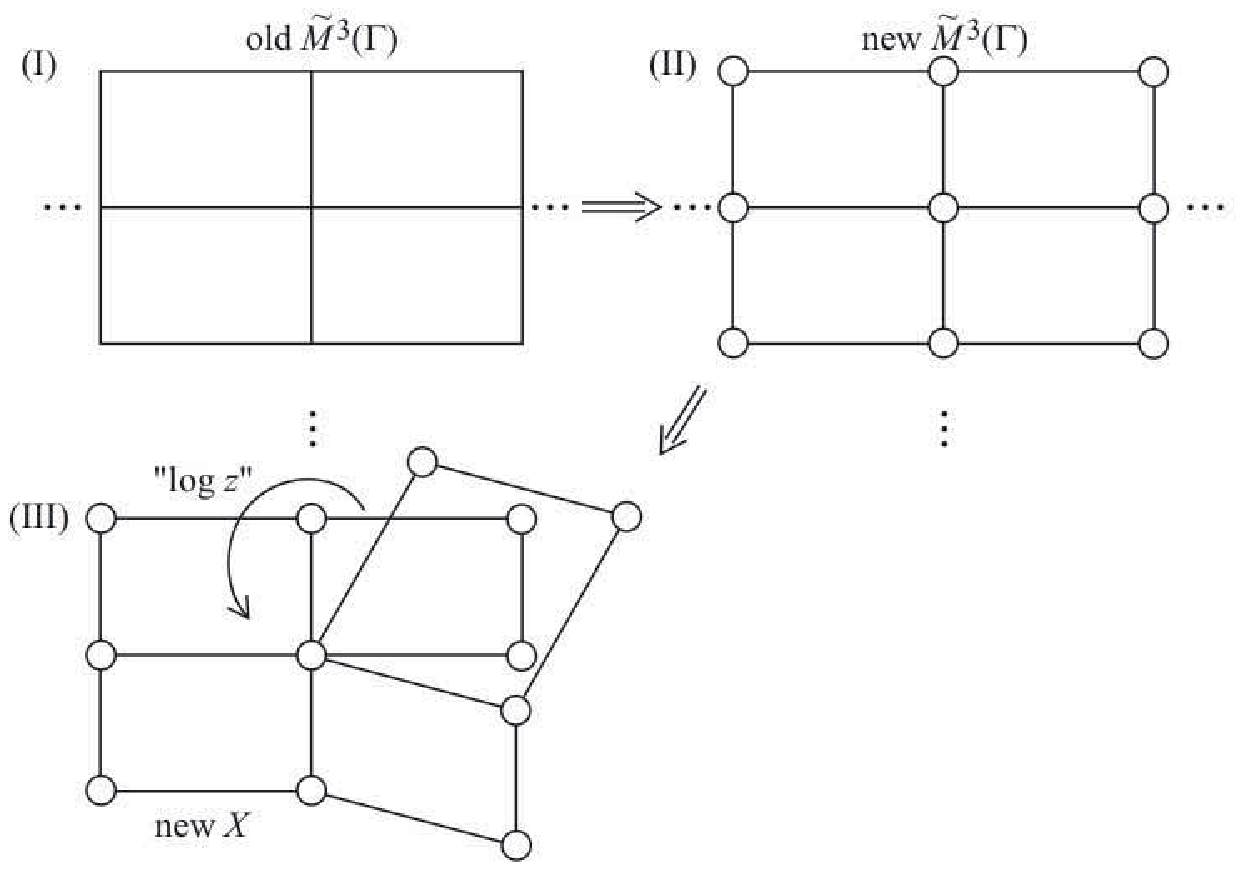}
$$
\centerline{Figure 2.1. This is a very schematical view of our redefinition old $\Rightarrow$ new.}
\begin{quote}
We have chosen here a toy-model where $M^3 (\Gamma)$ is replaced by $S^1 \times S^1$, with the old $\tilde M^3 (\Gamma)$ suggested in (I). When we go to the new context in (II), vertices become $0$-handles. The edges should become $1$-handles, but we have refrained from drawing this explicitly, so as to keep the figures simple. In (III) we are supposed to see singularities like in (\ref{eq1.15}). The pointwise $\log z$ singularity of the old $X$ is now an infinite chain of (\ref{eq1.15})-like singularities smeared in a $\log z$ pattern around a circle $S^1 \subset \delta (\mbox{$0$-handle})$. The notation ``$\log z$'' in (III) is supposed to suggest this.
\end{quote}

\bigskip

Our $X$ may also have a set of immortal singularities ${\rm Sing} \, (X) \subset X - {\rm Sing} \, (f)$ and it is exactly along ${\rm Sing} \, (f) + {\rm Sing} \, (X)$ that it fails to be a 3-manifold. It is along ${\rm Sing} (f)$ that it fails to be locally finite. Generically, the immortal singularities ${\rm Sing} (\tilde M^3 (\Gamma))$ are created by the zipping process and not mere images of the immortal ${\rm Sing} \, (X)$. Think of these latter ones as exceptional things, and with some extra work one could avoid them altogether. Anyway, immortal singularities are a novelty with respect to \cite{26}.

\smallskip

Also, in the next sections, the representation spaces will be completely devoid of immortal singularities.

\smallskip

The singularity $\sigma$ from (2.9.2) becomes now a finite linear chain of successive $0$-simplexes and $1$-simplexes. On the two extreme $0$-simplexes will occur now something like the ``$\log z$'' in figure~2.1-(III) which, outside them we have something like in the much more mundane figure 2.2-(II).

\smallskip

For the needs of the next sections, a bit more care has to go into the building of the singular handlebody $M^3 (\Gamma)$. To describe this, we reverse provisionally the transformation (2.11.1) and change $M^3 (\Gamma)$ back into a 2-dimensional cell-complex $[M^3 (\Gamma)]$, coming with $[\tilde M^3 (\Gamma)] = [M^3 (\Gamma)]^{\sim}$.

\bigskip

\noindent {\bf Lemma 2.5.} 1) {\it Without any loss of generality, i.e. without loosing any of our desirable features, our $M^3 (\Gamma)$ can be chosen such that $[M^3 (\Gamma)]$ is a {\ibf simplicial} complex.}

\smallskip

2) {\it At the same time, we can also make so that any given $2$-handle should  see at most {\ibf one} immortal singularity on its $\partial h^2$. Like before, this should involve another, distinct, $\partial h^2$.}

\smallskip

3) {\it All these properties are preserved by further subdivisions.}

$$
\includegraphics[width=13cm]{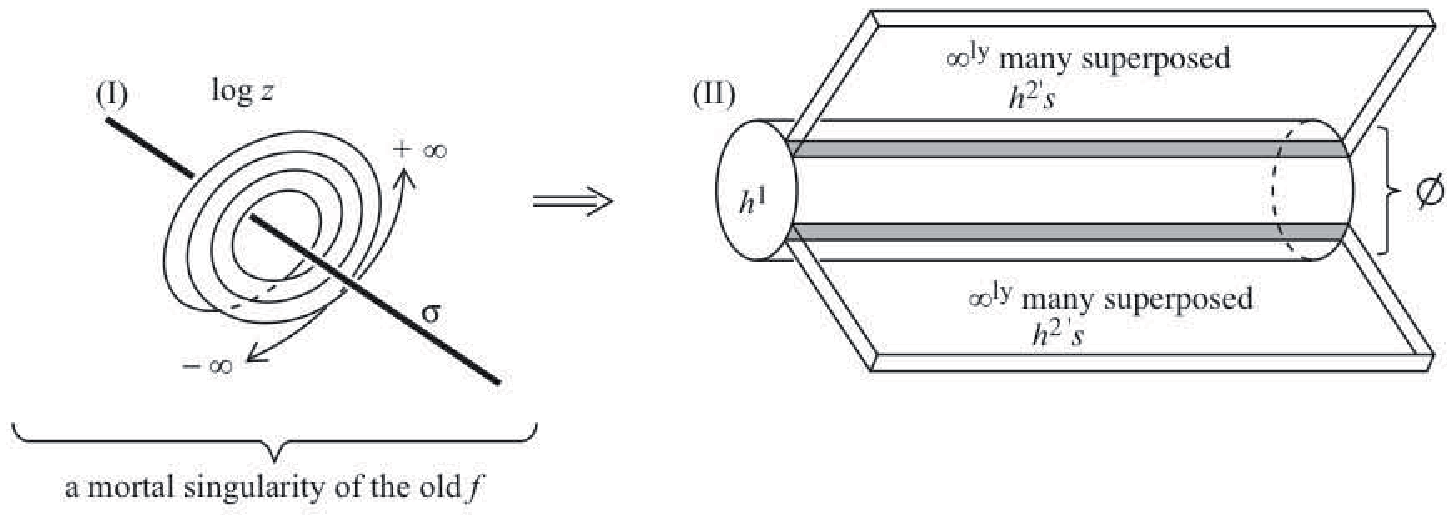}
$$
\centerline{Figure 2.2. We are here at the level of our $({\rm new}) \, X$.}
\begin{quote}
The $1$-handle $h^1$ corresponds to one of the small 1-simplices of $A^2 \mid \sigma$. The  hatched areas are mortal singularities $D^2 \subset {\rm Sing} \, (f)$.
\end{quote}

\bigskip

Here is how the easy lemma~2.5 is to be used. For convenience, denote $[\tilde M^3 (\Gamma)]$ by $K^2$ and let $K^1$ be its 1-skeleton. Consider then any arbitrary, non-degenerate simplicial map
\begin{equation}
\label{eq1.16}
S^1 \overset{\psi}{\longrightarrow} K^1 \, .
\end{equation}

\medskip

\noindent {\bf Lemma 2.6.} {\it The triangulation of $S^1$ which occurs in {\rm (\ref{eq1.16})} extends to a triangulation of $D^2$ such that there is now another simplicial {\ibf non-degenerate} map}
$$
D^2 \overset{\Psi}{\longrightarrow} K^2 \quad \mbox{with} \quad \Psi \mid S^1 = \psi \, .
$$

\medskip

The proof of lemma~2.6 uses the same argument as the proof of the ``simplicial lemma''~2.2 in \cite{26}. That kind of argument certainly has to use simplicial structures, hence the need for our lemma~2.5. I will only offer here some comments, without going into details.

\smallskip

The argument accompanying figure~2.7.3 in \cite{26} certainly needs requirement that two edges of $X^2$ with a common vertex should be joinable by a continuous path of 2-simplices. Without loss of generality, our $K^2$ verifies this condition (actually $[M^3 (\Gamma)]$ already does so). Next, it may be instructive to see what the argument boils down to when $X^2$ is reduced to a unique 2-simplex $\Delta$ the boundary of which is covered twice by (\ref{eq1.16}). Then, $D^2$ in the lemma is gotten out of four copies of $\Delta$. Start with a central $\Delta$, and then mirror it along its three sides. End of comment.

\section{Constructing equivariant locally-finite representations for $\tilde M^3 (\Gamma)$}\label{sec2}
\setcounter{equation}{0}

The present section will follow relatively closely \cite{26} and we will show how to extend that $3^{\rm d}$ part of \cite{26} which culminates with lemma~3.3, to $\tilde M^3 (\Gamma)$, which replaces now the smooth $\tilde M^3$ from \cite{26}. We will worry about the $2^{\rm d}$ context later on, in a subsequent paper. We use now again the formulae (\ref{eq1.13}) and (\ref{eq1.14}) from the last section, {\it i.e.} we write
\begin{equation}
\label{eq2.1}
\tilde M^3 (\Gamma) = \bigcup_{\lambda , i} h_i^{\lambda} \, , \quad X = \bigcup_{\lambda , i , \alpha} h_i^{\lambda} (\alpha) \, , \quad f (h_i^{\lambda} (\alpha)) = h_i^{\lambda} \, .
\end{equation}
We  may as well assume that $\Gamma$ operates already at the level of the indices $i$, so that
\begin{equation}
\label{eq2.2}
\mbox{For all $g \in \Gamma$, and $(i,\lambda)$ we have $g h_i^{\lambda} = h_{gi}^{\lambda}$ \, .}
\end{equation}
The $X \overset{f}{\longrightarrow} \tilde M^3 (\Gamma)$ will be replaced now by a new $3^{\rm d}$ representation
\begin{equation}
\label{eq2.3}
Y \overset{G}{\longrightarrow} \tilde M^3 (\Gamma)
\end{equation}
having the properties 1) (local finiteness) and 2) (equivariance) from our REPRESENTATION THEOREM. We will worry about the bounded zipping length only later on in the next section. Our $Y$ which is a Thurston $3^{\rm d}$ train-track is a handlebody type union of {\ibf bicollared handles} $H_i^{\lambda} (\gamma)$, where $\gamma$ belongs to a countable family of indices, a priori $(\lambda , i)$-dependent. Bicollared handles were defined in \cite{26} and, according to this definition, a bicollared handle $H$ is, purely topologically speaking, a bona fide handle $\hat H$ from which the lateral surface (which we will denote by $\delta \hat H$, or simply by $\delta H$) has been {\ibf deleted}. But bicollared handles have more structure, namely a filtration by bona fide handles with the same index as $H$, $H = \overset{\infty}{\underset{n=1}{\bigcup}} \, H_n$ (see (2.15) in \cite{26}). This endows $H$ with two collars, each with countably many layer (or ``levels'') and incoming collar parallel to the attaching zone $\partial H^{\lambda}$ of $H^{\lambda}$ and an outgoing collar, parallel to the lateral surface $\delta H^{\lambda}$. These collars, which are not disjoined, can be visualized in the figures 2.1, 2.3 from \cite{26}. The general idea is that, when $H^{\lambda}$ is attached to $H^{\lambda - 1}$ along (a piece of) $\partial H^{\lambda}$, making use of the respective outgoing collar of $H^{\lambda - 1}$ and incoming collar of $H^{\lambda}$, then the two outgoing collars, of $H^{\lambda - 1}$ and $H^{\lambda}$ {\ibf combine} into a unique outgoing collar for $H^{\lambda - 1} \cup H^{\lambda}$. It is that part of the outgoing collar of $H^{\lambda - 1}$ which was not used for attaching $H^{\lambda}$ which occurs here. The physical glueing of $H^{\lambda}$ to $H^{\lambda - 1}$ occurs actually along a PROPERLY embedded codimension one hypersurface $\partial H^{\lambda} \cap H^{\lambda - 1} \subset H^{\lambda - 1}$.

\smallskip

It is along the newly created collar above that $H^{\lambda + 1}$ will be attached. Moreover, when a handle is to be attached along some outgoing collar this always happens at some specific {\ibf level} $\ell$, see here \cite{26}. We will want to control these levels.

\smallskip

Assume, for instance, that we have a necklace of bicollared $0$-handles and $1$-handles, where the bicollared $H_i^1$ is attached at level  $\ell'_i$ to its left and $\ell''_i$ to its right. As explained in \cite{26} (see the pages 28 and 29), the following {\ibf ``frustration number''} ($\approx$ holonomy, in a discrete version)
$$
K = \sum_i (\ell''_i - \ell'_{i+1})
$$
is the obstruction for the necklace to have a good outgoing collar to which a bicollared $2$-handle can be attached. This obstruction is very easily dealt with by observing first that $K$ may be re-written as $\underset{i}{\sum} \, (\ell'_i - \ell''_i)$ and then making sure that when $H_i^1$ is attached, we always fix the levels so that $\ell'_i = \ell''_i$.

\smallskip

As far as (\ref{eq2.3}) is concerned, all this is $Y$-story. Next, the map $G$ is always supposed to be such that each $G \mid H_i^{\lambda} (\gamma)$ extends continuously to a larger embedding $G \mid \hat H_i^{\lambda} (\gamma)$, where $\hat H_i^{\lambda} (\gamma) = H_i^{\lambda} (\gamma) \cup \delta H_i^{\lambda} (\gamma)$. Inside $\tilde M^3 (\Gamma)$, for each $h_i^{\lambda}$, the $GH_i^{\lambda} (\gamma)$ occupies, {\it roughly}, the position $h_i^{\lambda}$ while, again for each $h_i^{\lambda}$, we always have the {\ibf strict} equality $G \delta \hat H_i^{\lambda} (\gamma) = \delta h_i^{\lambda}$, for all $\lambda$'s.

\smallskip

Using the technology from \cite{26} we can get now the following lemma, as well as the two complements which follow.

\bigskip

\noindent {\bf Lemma 3.1.} {\it We can perform the construction of {\rm (\ref{eq2.3})} so that}

\smallskip

1) {\it $Y$ is GSC, $\Psi (G) = \Phi (G)$ (i.e. $G$ is zippable) and also
$$
\overline{{\rm Im} \, G} = \tilde M^3 (\Gamma) \, .
$$
Of course, so far this only expresses the fact that} (\ref{eq2.3}) {\it is a representation, but then we also have the next items.}

\smallskip

2) {\it A FIRST FINITENESS CONDITION (at the source). The complex $Y$ is {\ibf locally finite}.}

\smallskip

3) {\it There is a {\ibf free action} $\Gamma \times Y \to Y$, for which the map $G$ is {\ibf equivariant}, i.e.
$$
G(gx) = g \, G(x) \, , \quad \forall \, x \in Y \, , \ g \in \Gamma \, .
$$
Moreover, with the same action of $\Gamma$ on the indices $i$ like in} (\ref{eq2.2}), {\it the action of $\Gamma$} on the $0$-skeleton {\it of $Y$ takes the following form, similar to} (\ref{eq2.2}), {\it namely
$$
g H_i^0 (\gamma) = H_{gi}^0 (\gamma) \, , \quad \forall \, \gamma \, .
$$
By now we really have an equivariant, locally-finite representation for our $\tilde M^3 (\Gamma)$ ($\approx \Gamma$), like it is announced in the title of the present section.}

\smallskip

4) {\it We have $GH_i^0 (\gamma) = {\rm int} \, h_i^0$ and also $G(\delta H_i^0 (\gamma)) = \delta h_i^0$, as already said above. This makes the $GH_i^0 (\gamma)$ be $\gamma$-independent, among other things.}

\bigskip

Retain that the $\delta H_i^{\lambda} (\gamma)$ exists only ideally, as far as the $Y$ is concerned; it lies at infinity. But we certainly can introduce the following stratified surface, which is PROPERLY embedded inside $\tilde M^3 (\Gamma)$, namely
\begin{equation}
\label{eq2.4}
\Sigma_1 (\infty) \underset{\rm def}{=} \bigcup_{i,\lambda , \gamma} G (\delta H_i^{\lambda} (\gamma)) = \bigcup_{i,\lambda} \delta h_i^{\lambda} \subset \tilde M^3 (\Gamma) \, ,
\end{equation}
which comes with its useful restriction
\begin{equation}
\label{eq2.5}
\Sigma_2 (\infty) \underset{\rm def}{=} GY \cap \Sigma_1 (\infty) = \bigcup \ \{ \mbox{common faces $h_i^{\lambda} \cap \lambda_j^{\mu}$,} 
\end{equation}
$$
\mbox{for $0 \leq \lambda < \mu \leq 2$}\} = \bigcup \ \{\mbox{interiors of the attaching zones} \ \partial h^1 , \partial h^2\}.
$$
In the lemma which follows next, the $\varepsilon$-skeleton of $Y$ is denoted $Y^{(\varepsilon)}$.

\bigskip

\noindent {\bf Lemma 3.2.} (FIRST COMPLEMENT TO LEMMA 3.1.) 1) {\it The $Y^{(\varepsilon)}$ contains a canonical outgoing collar such that each $H_i^{\lambda} (\gamma)$ is attached to $Y^{(\lambda - 1)}$ in a collar-respecting manner at some level $k (i,\gamma) \in Z_+$ which is such that}
\begin{equation}
\label{eq2.6}
\lim_{i+\gamma = \infty} k(i,\gamma) = \infty \, .
\end{equation}

2) {\it For each $Y^{(\varepsilon)}$ we introduce its ideal boundary, living at infinity, call it
$$
\delta Y^{(\varepsilon)} = \bigcup_{\overbrace{i,\gamma,\lambda \leq \varepsilon}} \delta H_i^{\lambda} (\gamma) \quad \mbox{and also} \quad \hat Y^{(\varepsilon)} = Y^{(\varepsilon)} \cup \delta Y^{(\varepsilon)} \, .
$$
Eeach $H_i^{\lambda} (Y)$ is attached to $Y^{(\lambda - 1)}$ via $\partial H_i^{\lambda} (\gamma)$ in a bicollared manner, and as a consequence of} (\ref{eq2.6}) {\it inside $\hat Y^{(\lambda - 1)}$ we will find that}
$$
\lim_{n+m=\infty} \partial H_m^{\lambda} (\gamma_n) \subset \delta Y^{(\lambda - 1)} \, , \eqno (3.6.1)
$$
{\it which implies the FIRST FINITENESS CONDITION from lemma}~3.1.

\smallskip

3) {\it Inside $\tilde M^3 (\Gamma)$ we also have that}
$$
\lim_{n+m+k=\infty} G \delta H_{n,m}^{\lambda} (\gamma_k) \subset G \delta Y^{(\lambda)} = \tilde M^3 (\Gamma)^{(\lambda)} \subset \Sigma_1 (\infty) \, . \eqno (3.6.2)
$$

\medskip

\noindent A COMMENT. In point 4) of the 2-dimensional representation theorem which was stated in the introduction (but which will only be proved in a subsequent paper) we have introduced the so-called SECOND FINITENESS CONDITION. Eventually, it will be the (3.6.2) which will force this condition.

\bigskip

\noindent {\bf Lemma 3.3.} (SECOND COMPLEMENT TO LEMMA 3.1.) 1) {\it There are PROPER individual embeddings $\partial H_i^{\lambda} (\gamma) \subset Y^{(\lambda - 1)}$ and the global map
\begin{equation}
\label{eq2.7}
\sum_{0 < \lambda \leq 2,i,\gamma} \partial H_i^{\lambda} (\gamma) \overset{j}{\longrightarrow} Y
\end{equation}
is also PROPER. As a novelty with respect to} \cite{26}, {\it the map $j$ above fails now to be injective. At each immortal singularity $S$ of $X$ (if such exist), call it $S \subset \partial h_k^0 (\alpha_1) \cap \partial h_i^2 (\alpha_2) \cap \partial h_j^2 (\alpha_3)$, we have transversal contacts $\partial H_i^2 (\alpha_2) \pitchfork \partial H_j^2 (\alpha_3) \subset H_k^0 (\alpha_1)$.}

\smallskip

2) {\it We have ${\rm Im} \, j = {\rm Sing} \, (G)$ ($=$ mortal singularities of $G$) and, at the same time, ${\rm Im} \, j$ is the set of non-manifold points of the traintrack $Y$.}

\smallskip

3) {\it There are {\ibf no} immortal singularities for $Y$, all the immortal singularities for $GY$ are created by the zipping.}

\bigskip

We will describe now the GEOMETRY OF $(G,Y)$ in the neighbourhood of an immortal singularity downstairs
$$
S \subset \delta h_k^0 \cap \partial h_i^2 \cap \partial h_j^2 \subset \tilde M^3 (\Gamma) \, .
$$
When we fix the indices $k,i,j$ then, at the level of $Y$ we will find an infinity of triplets $H_k^0 (\gamma_n)$, $H_i^2 (\gamma'_n)$, $H_j^2 (\gamma''_n)$, $n \to \infty$. Here each of the $G \, \underset{n}{\sum} \, H_k^0 (\gamma_n) \cup H_i^2 (\gamma'_n)$, $G \, \underset{n}{\sum} \, H_k^0 (\gamma_n) \cup H_j^2 (\gamma''_n)$ generates a figure analogous to 2.2 in \cite{26}, living inside the common $GH_k^0 (\gamma_n) = {\rm int} \, h_k^0$. What we are discussing now is the interaction of these two figures.

\smallskip

For any given pair $n,m$, among the two $\partial H_i^2 (\gamma'_n)$, $\partial H_j^2 (\gamma''_m)$ one is ``low'', {\it i.e.} pushed deeper inside ${\rm int} \, h_k^0$, the other one is ``high'', {\it i.e.} living more shallowly close to $\delta h_k^0$.

\smallskip

If we think, loosely of each $\partial H^2$ as being a copy of $S^1 \times [0,1]$ then the intersection $G \partial H_i^2 (\gamma'_n) \cap G \partial H_j^2 (\gamma''_m)$ with $\gamma'_n$ low and $\gamma''_n$ high, consists of two transversal intersection arcs which are such that, from the viewpoint of the lower $\gamma'_n$ they are close to the boundary $S^1 \times \{ 0,\varepsilon \}$, while from the viewpoint of the higher $\gamma''_m$ they are generators $p \times [0,1]$ and $q \times [0,1]$.

\smallskip

Given $\gamma'_n$, for almost all $\gamma''_m$ the $\gamma'_n$ is low and the $\gamma''_m$ is high, and a similar thing is true when we switch $\gamma'$ and $\gamma''$.

\bigskip

\noindent SKETCH OF PROOF FOR THE LEMMA 2.1 AND ITS COMPLEMENTS. We will follow here rather closely \cite{26}. Together with the transformation $\{ (X,f)$ $\mbox{from (3.1)} \Longrightarrow \{ (Y,G) \ \mbox{from (3.3)} \}$, to be described now, will come two successive increases of the family of indices
\begin{equation}
\label{eq2.8}
\underbrace{\{ \alpha \}}_{{\rm like \ in \ (3.1)}} \subsetneqq \{ \beta \} \subsetneqq \underbrace{\{ \gamma \}}_{{\rm like \ in \ (3.3)}} \, .
\end{equation}
The first step, is to change each $h_i^{\lambda} (\alpha) \subset X$ into a bicollared handle $H_i^{\lambda} (\alpha)$. The $H_i^{\lambda} (\alpha)$'s are put together into a provisional $Y = Y (\alpha)$ following roughly, but with some appropriate perturbations, the recipee via which the $h_i^{\lambda} (\alpha)$'s themselves have been put together at the level of $X$ (3.1). 

\smallskip

The index ``$\alpha$'' is here generic, standing for ``first step'', and it is not to be mixed up with the ``$\alpha$'' which occurs inside $H_i^{\lambda} (\alpha)$. But then, when this index is used for $Y$ rather than for an individual $H^{\lambda}$ we may safely write again ``$\alpha$''.

\smallskip

For the time being we will concentrate on $\lambda = 0$, where we start by {\ibf breaking the $i$-dependence} of the $\alpha$'s, chosing for each $i$ a fixed isomorphism $\{ \alpha \} \approx Z_+$. For $\lambda = 0$, the index ``$\alpha$'' has now a universal meaning; this will be essential for the (\ref{eq2.10}) below. We consider a provisional, restricted $0$-skeleton for $Y$,
\begin{equation}
\label{eq2.9}
Y^{(0)} (\alpha) = \sum_{i,\alpha} H_i^0 (\alpha) \, , \quad \mbox{with} \ G \hat H_i^0 (\alpha) = h_i^0 (\alpha) \, .
\end{equation}
Next, we will {\ibf force} a free action $\Gamma \times Y^{(0)} (\alpha) \to Y^{(0)} (\alpha)$, by
\begin{equation}
\label{eq2.10}
g H_i^0 (\alpha) = H_{gi}^0 (\alpha) \, \quad \forall \, g \in \Gamma \, , \ \mbox{with the {\ibf same} index $\alpha$ on both sides.}
\end{equation}
Whatever additional requirements there will be for the collaring, we will also always have
\begin{equation}
\label{eq2.11}
gH_{n,i}^0 (\alpha) = H_{n,gi}^0 (\alpha) \, , \quad \forall \, n \, .
\end{equation}
With (\ref{eq2.2}), (\ref{eq2.10}), (\ref{eq2.11}) comes easily a $\Gamma$-equivariant $G \mid Y^{(0)} (\alpha)$, and the following condition will be imposed on our equivariant collaring too.

\smallskip

The following set accumulates, exactly, on $\delta h_i^0$
\begin{equation}
\label{eq2.12}
\underbrace{\sum_{n,\alpha} G (\delta H_{n,i}^0 (\alpha))}_{\mbox{this is the same thing as} \atop \mbox{$\underset{g,n,\alpha}{\sum} g^{-1} G (\delta H_{n,gi}^0 (\alpha))$}} \!\!\!\!\!\!\!\!\!\!\!\!\subset h_i^0 \subset \tilde M^3 (\Gamma) \, .
\end{equation}
We turn now to the $1$-skeleton. At the level of $X$ (\ref{eq2.1}), we have things like
\begin{equation}
\label{eq2.13}
\partial h_j^1 (\alpha) = h_{j_0}^0 (\alpha_0) - h_{j_1}^0 (\alpha_1) \, ,
\end{equation}
where each pair $(j_{\varepsilon} , \alpha_{\varepsilon})$ depends of $(j,\alpha)$. With this, at the level of our provisional, restricted $1$-skeleton
$$
Y^{(1)} (\alpha) = Y^{(0)} (\alpha) \cup \sum_{i,\alpha} H_i^1 (\alpha) \subset Y \ (\mbox{from (3.3)}) \, ,
$$
we will attach $H_i^1 (\alpha)$ to each of the two $H_{j_{\varepsilon}}^0 (\alpha_{\varepsilon})$ at some respective level $k (j_{\varepsilon} , \alpha_{\varepsilon})$, for which we will impose the condition
\begin{equation}
\label{eq2.14}
k (j_0 , \alpha_0) = k (j_1 , \alpha_1) \, ,
\end{equation}
which will make sure that for the frustration numbers we have $K=0$. This will be a standard precautionary measure, from now on. The (\ref{eq2.14}) is to be added to the conditions (\ref{eq2.6}), of course.

\smallskip

But, of course, there is no reason for our $\underset{i,\alpha}{\sum} \ H_i^1 (\alpha)$ to be equivariant, contrary to what happens for the $\underset{i}{\sum} \ h_i^1 \subset \tilde M^3 (\Gamma)$. We will {\ibf force} now such an equivariance by the following kind of averaging procedure. We start by extending the family $\{ \alpha \}$ into a larger family of indices $\{ \beta \} \supsetneqq \{ \alpha \}$ such that, for the time being at a purely abstract level, the following {\ibf saturation formula} should be satisfied
\begin{equation}
\label{eq2.15}
\{\mbox{the set of all $H_h^1 (\beta)$, $\forall \, h$ and $\beta\} = \{$the set of all $gH_j^1 (\alpha)$, $\forall \, g,j$ and $\alpha\}$.}
\end{equation}
Here $(h,\beta) = (h,\beta) [g,j,\alpha]$ and, at the purely abstract level of the present discussion, we will also impose things like
$$
g_1 \, H_h^1 (\beta) = g_1 \, g H_j^1 (\alpha) = H_{h(g_1 g,j,\alpha)}^1 (\beta (g_1 \, g , j,\alpha)) \ \mbox{and} \ (h,\beta) (1 \in \Gamma , j,\alpha) = (j,\alpha) \, .
$$

The next step is now to give flesh and bone to the abstract equality of sets (\ref{eq2.15}) by requiring that, for each $(g,j,\alpha)$, the corresponding $H_h^1 (\beta)$ (from (\ref{eq2.15})) should be attached to the two $0$-handles $g H_{j_{\varepsilon}}^0 (\alpha_{\varepsilon})$ (see (\ref{eq2.13})), and this at the {\ibf same} level $k (j_{\varepsilon} , \alpha_{\varepsilon})$ (\ref{eq2.14}) on both sides. This way we have created an object $Y^{(1)} (\beta) \supset Y^{(0)} (\beta) \underset{\rm def}{=} Y^{(0)} (\alpha)$, which is endowed with a free action of $\Gamma$, which extends the preexisting action on $Y^{(0)} (\alpha)$. We will also ask that the saturation should work at the level of the $H_{n,i}^1 (\alpha)$'s too. Moreover we can arrange things so that the following implication should hold
$$
\{ g H_i^1 (\beta) = H_{i_1}^1 (\beta_1)\} \Longrightarrow \{ g H_{n,i}^1 (\beta) = H_{n,i_1}^1 (\beta_1) \ \mbox{for all $n$'s}\} .
$$
By now we can easily get an equivariant map too, extending $G \mid Y^{(0)} (\alpha)$
\begin{equation}
\label{eq2.16}
Y^1 (\beta) \overset{G}{\longrightarrow} (\tilde M^3 (\Gamma))^{(1)} = \sum_{\lambda \leq 1,i} h_i^{\lambda} \, .
\end{equation}
The map (\ref{eq2.16}) is not completely rigidly predetermined, inside $\tilde M^3 (\Gamma)$ the $GH_j^1 (\beta)$ occupies only {\it roughly} the position $h_j^1$, but then we will also insist on strict equalities $G\delta H_j^1 (\beta) = \delta h_j^1$. As a side remark, notice that given $j,\beta , g$, there is some $\beta' = \beta' (j,\beta ,g)$ with the equality
$$
g H_j^1 (\beta) = H_{gj}^1 (\beta') \, , \ \mbox{corresponding to} \ gh_j^1 = h_{gj}^1 \ \mbox{(see (3.2))}.
$$
Finally, we can also impose on (\ref{eq2.16}) the condition that the doubly infinite collection of annuli
$$
\sum_{g,n,\beta} g^{-1} G (\delta H_{n,gi}^1 (\beta)) = \sum_{n,\beta} G (\delta H_{n,i}^1 (\beta)) \subset (\tilde M^3 (\Gamma))^{(1)}
$$
accumulates exactly on $\delta h_i^1$, just like in (\ref{eq2.12}).

\smallskip

Once all the {\it frustration numbers are zero}, we can next proceed similarly for the $2$-skeleton and thereby build a completely equivariant $Y(\beta)$, endowed with an equivariant map into $\tilde M^3 (\Gamma)$. At this level, all the requirements in our lemma~3.1 (and its complements) are fulfilled, {\ibf except} the geometric sample connectivity which, until further notice has been violated. To be very precise, the GSC, which was there at the level of $X$ (\ref{eq2.1}), got lost already when at the level of the saturation formula  (\ref{eq2.15}), additional $1$-handles were thrown in. There is no reason, of course, why the $\{ H_j^2 (\beta) \}$ should cancell them. Before going on with this, here are some more details concerning the construction $Y^{(1)} (\beta) \Longrightarrow Y(\beta)$, which yields the full 2-skeleton. We start with an abstract saturation formula, analogous to (\ref{eq2.15})
$$
\{ H_h^2 (\beta) \} = \{ g H_i^2 (\alpha) \} \, , \ \forall \, h, \beta , g ,i,\alpha \, ,
$$
producing 2-handles which can afterwards be happily attached to $Y^{(1)} (\beta)$.

\smallskip

So, we have by now a map which is equivariant, has a locally finite source, but with a source which fails to be GSC, reason which prevents it to be a representation,
$$
Y(\beta) \overset{G}{\longrightarrow} \tilde M^3 (\Gamma) \, .
$$
We want to restore the GSC, without loosing any of the other desirable features already gained. Proceeding like in \cite{26}, inside the 1-skeleton $Y^{(1)} (\beta)$ we chose a PROPER, equivariant family of simple closed curves $\{ \gamma_j (\beta) \}$, in cancelling position with the family of 1-handles $\{ H_j^1 (\beta) \}$. With this, we extend $\{ G \gamma_j (\beta) \}$ to an equivariant family of mapped discs
\begin{equation}
\label{eq2.18}
\sum_{j,\beta} D^2 (j,\beta) \overset{f}{\longrightarrow} \tilde M^3 (\Gamma) \, ,
\end{equation}
and the first idea would be now to replace $Y(\beta)$ by the following equivariant map
\begin{equation}
\label{eq2.19}
Y(\beta) \cup \sum_{j,\beta} D^2 (j,\beta) \overset{G \cup f}{-\!\!\!-\!\!\!-\!\!\!-\!\!\!-\!\!\!\longrightarrow} \tilde M^3 (\Gamma) \, ,
\end{equation}
the source of which is both locally finite and GSC. We may also assume that $Y(\beta) \overset{G}{\longrightarrow} \tilde M^3 (\Gamma)$ is essentially surjective, from which it is easy to get $\Psi (G \cup f) = \Phi (G \cup f)$, {\it i.e.} a zippable $G \cup f$. In other words (\ref{eq2.19}) is already a representation, even an equivariant, locally finite one. BUT there is still a trouble with (\ref{eq2.19}), namely that we do not control the accumulation pattern of the family $\{ f D^2 \}$ which, a priori, could be as bad as {\ibf transversally cantorian}. In particular, the SECOND FINITENESS CONDITION (\ref{eq0.6}), which will be essential in the subsequent paper of this series is, generically speaking, violated by (\ref{eq2.19}). [More precisely, as one can see in \cite{26} and \cite{27}, once representations go $2$-dimensional, the accumulation pattern of $M_2 (f)$, instead of being, transversally speaking, with finitely many limit point like (\ref{eq0.6}) wants it to be, tends to become cantorian.]

\smallskip

The cure for this disease is the following {\ibf decantorianization process}, which is discussed at more length in \cite{26}. Here are the successive steps 

\bigskip

\noindent (3.18.1) We start by applying lemma~2.6 to $\{ G \gamma_j (\beta) \}$, and this replaces (\ref{eq2.18}) by a simplicial {\ibf non degenerate} map
$$
\sum_{j,\beta} D^2 (j,\beta) \overset{f}{\longrightarrow} [M^3 (\Gamma)] \quad \mbox{(see section II)}.
$$
We also have the following features

\bigskip

\noindent (3.18.2) \quad The $\underset{j,\beta}{\sum} \ D^2 (j,\beta)$ is endowed with an equivariant triangulation, involving three infinite families of simplices $\{ s^0 \}$, $\{ s^1 \}$, $\{ s^2 \}$ of dimensions 0, 1, and 2, respectively.

\bigskip

\noindent (3.18.3) \quad For any simplex $k = s^{\lambda}$ not on the boundary and not corresponding hence to some already existing bicollared handle $H_i^{\lambda} (\beta)$, we introduce a new bicollared family $H_k^{\lambda} (\gamma)$. Next, we use $Y(\beta) \cup \underset{j,\beta}{\sum} \ D^2 (j,\beta)$ as a pattern for attaching, successively, these new bicollared $H_k^{\lambda} (\gamma)$ to $Y(\beta)$. This leads to a $3^{\rm d}$ representation of $\Gamma$,
\begin{equation}
\label{eq2.20}
Y = Y (\gamma) = \bigcup_{\lambda , i , \gamma} H_i^{\lambda} (\gamma) \overset{G}{\longrightarrow} \tilde M^3 (\Gamma) \, .
\end{equation}
We follow here the already established patterns of controlling the levels when bicollared handles are attached, so as to have local finiteness (FIRST FINITENESS CONDITION) and, at the same time the SECOND FINITENESS CONDITION. The equivariance of (\ref{eq2.20}) will descend directly from the equivariance of (\ref{eq2.19}), and similarly the condition $\Psi = \Phi$ too. With this, our lemmas~3.1, 3.2, 3.3 can be established.

\bigskip

\noindent FINAL REMARKS. A) In view of (\ref{eq2.8}) we have here also inclusions at the level of the handles
$$
\{ H_i^{\lambda} (\alpha) \} \subsetneqq \{ H_i^{\lambda} (\beta) \} \subsetneqq \{ H_i^{\lambda} (\gamma) \} \, ,
$$
hence $Y(\alpha) \subset Y(\beta) \subset Y(\gamma)$, with compatible maps into $\tilde M^3 (\Gamma)$.

\smallskip

\noindent B) Here is a very fast way to understand how one decantorianizes via a process of {\ibf thickening}. Inside the interval $[a,b]$ far from the border, consider a sequence $x_i \in (a,b)$, $i=1,2,\ldots$ accumulating on a Cantor set $C \subset (a,b) - \underset{i}{\sum} \ \{ x_i \}$.

\smallskip

In order to get rid of $C$, we start by consider two sequences $a_n , b_n \in (a,b) - C - \{ x_i \}$, such that $\lim a_n = a$, $\lim b_n = b$. If one thickens each $x_n$ int $[a_n , b_n] \ni x_n$ then one has managed to replace the Cantorian accumulation of the set $\{ x_n \}$ by the very tame accumulation of the arcs $\{ [a_n , b_n ] \}$.

\smallskip

Of course, what we are worrying about is the transversal accumulation to an infinite family of $2^{\rm d}$ sheets in $3^{\rm d}$. There, one decantorianizes by thickening the $2^{\rm d}$ sheets into $3^{\rm d}$ boxes (or $3^{\rm d}$ handles of index $\lambda = 2$). With a bit of care, the accumulation pattern of the lateral surfaces of the handles is now tame.

\section{Forcing a uniformaly bounded zipping length}\label{sec3}
\setcounter{equation}{0}

The last section has provided us with the 3-dimensional representation (\ref{eq2.3}), {\it i.e.} with 
$$
Y \overset{G}{\longrightarrow} \tilde M^3 (G)
$$
which had both a locally finite representation space $Y$ and was equivariant. We turn now to the issue of the zipping length and, starting from (\ref{eq2.3}) as a first stage, what we will achieve now will be the following

\bigskip

\noindent {\bf Lemma 4.1.} {\it There is another, much larger $3$-dimensional representation, extending the} (\ref{eq2.3}) {\it above,
\begin{equation}
\label{eq3.1}
\xymatrix{
Y(\infty) \ar[rr]^{\!\!\!\!\!\!\!\!\!\!\ g(\infty)}  && \ \tilde M^3 (\Gamma) \, , 
} 
\end{equation}
having all the desirable features of} (\ref{eq2.3}), {\it as expressed by lemma~{\rm 3.1} and its two complements} 3.2, 3.3 {\it and moreover such that the following three things should happen too.}

\bigskip

\noindent (4.2) \quad {\it There exists a {\ibf uniform bound} $M > 0$ such that for any $(x,y) \in M^2 (g(\infty))$ we have
$$
\inf_{\lambda} \Vert \lambda (x,y) \Vert < M \, ,
$$
where $\lambda$ runs through all the zipping paths for $(x,y)$.}

\bigskip

\noindent (4.3) \quad {\it Whenever $(x,y) \in M^2 (g(\infty)) \mid Y(\infty)^{(1)}$, then we can find a zipping path of controlled length $\Vert \lambda (x,y) \Vert < M$, which is confined inside the $1$-skeleton $Y(\infty)^{(1)}$.}

\bigskip

\noindent (4.4) \quad {\it When $\Psi (g(\infty)) = \Phi (g(\infty))$ is being implemented by a complete zipping of $g(\infty)$, then one can start by zipping completely the $1$-skeleton of $Y(\infty)$ before doing anything about the $2$-skeleton.}

\bigskip

\noindent PROOF OF THE LEMMA 4.1. We start with some comments. The bounded zipping length is really a $1$-skeleton issue. Also, it suffices to check it for double points $(x,y) \in M^2 (g(\infty)) \mid Y (\infty)^{(0)}$ with a controlled $\lambda (x,y)$ confined inside $Y(\infty)^{(1)}$.

\smallskip

We will use the notation $(Y^{(0)} , G_0) \underset{\rm def}{=} (Y,G)$ (\ref{eq2.3}) and, with this, our proof will involve an infinite sequence of $3$-dimensional representation $(Y (n) , G_n)$ for $n = 0,1,2,\ldots$, converging to the desired $(Y(\infty) , g(\infty))$. We will find
$$
Y(0) \subset Y(1) \subset Y(2) \subset \ldots , \quad Y (\infty) = \bigcup_{0}^{\infty} Y(n)
$$
where the $Y(0) = Y$ is already GSC and where the construction will be such that the $2$-handles in $Y(n) - Y(n-1)$ will be in cancelling position with the $1$-handles in $Y(n) - Y(n-1)$, something which will require the technology of the previous section. Anyway, it will follow from these things that $Y(\infty) \in {\rm GSC}$.

\smallskip

The construction starts by picking up some fundamental domain $\Delta \subset \tilde M^3 (\Gamma)$ of size $\Vert \Delta \Vert$. Let us say that this $\Delta$ consists of a finite connected system of handles $h^0 , h^1 , h^2$. No other handles will be explicitly considered, until further notice, in the discussion which follows next. One should also think, alternatively, of our $\Delta$ as being part of $M^3 (\Gamma)$.

\smallskip

We will start by establishing abstract, intercoherent isomorphisms between the various, a priori $\ell$-dependent systems of indices $\gamma$, for the various infinite families $\{ H_{\ell}^0 (\gamma) \}$ which correspond to the finitely many $h_{\ell}^0 \subset \Delta$. Downstairs, in $\Delta$, we also fix a triple $h_1^0 , h_2^0 , h_i^1$ such that $\partial h_i^1 = h_2^0 - h_1^0$. Upstairs, some specific $H_1^0 (\gamma_0)$ living over $h_1^0$ will be fixed too.

\smallskip

With all this, we will throw now into the game a system of bicollared handles $H_j^1 (\delta)$, living over the various $h_j^1 \subset \Delta$, many enough so as to realize the following conditions

\bigskip

\noindent (4.5) \quad For any given index $\gamma$ there is an index $\delta = \delta (\gamma)$, such that $\partial H_i^1 (\delta) = H_1^0 (\gamma_0) - H_2^0 (\gamma)$, with $H_i^1 (\delta)$ living over $h_i^1$.

\bigskip

\noindent (4.6) \quad Downstairs, for each pair $h_n^0 , h_m^0 \subset \Delta$ we fix a path joining them $\mu = \mu (h_n^0 , h_m^0) \subset \tilde M^3 (\Gamma)$, with $\Vert \mu \Vert \leq \Delta$. It is assumed here that $\mu (h_1^0 , h_2^0) = h_i^1$. With this, for any index $\gamma$ and each pair $h_n^0 , h_m^0 \subset \Delta$ it is assumed that enough things should have been added upstairs, s.t. there should be a continuous path $\lambda = \lambda (h_n^0 , h_m^0 , \gamma)$ covering exactly the $\mu (h_n^0 , h_m^0)$ and such that, moreover,
$$
\partial \lambda (h_n^0 , h_m^0 , \gamma) = H_n^0 (\gamma) - H_m^0 (\gamma) \, , \ \forall \, \gamma \, .
$$
Here $\lambda (h_1^0 , h_2^0 , \gamma_0) = H_i^1 (\delta (\gamma_0))$. From here on we will imitate the process $Y(\alpha) \subset Y(\beta) \subset Y(\gamma) = \{ Y ,$ {\it i.e.} our present $Y(0)\}$, meaning the following succession of steps. The only addition, so far, has been the system of $1$-handles $\{ H_1^0 (\delta)\}$, localized over $\Delta$. The next step will be to saturate it like in (\ref{eq2.15}) so as to achieve equivariance. The next two steps will add more handles of index $\lambda = 0$, $\lambda = 1$ and $\lambda = 2$ so as to restore GSC and then to decantorianize. Each of these two steps has equivariance built in and also, by proceeding like in the last section, local finiteness is conserved. The feature $\Psi = \Phi$ is directly inherited from $Y(0)$, whose image is already almost everything.

\smallskip

The end result of this process is a bigger $3$-dimensional representation
$$
Y \underset{\rm def}{=} Y(0) \subset Y(1) \overset{G_1}{-\!\!\!-\!\!\!\longrightarrow} \tilde M^3 (\Gamma) \, , \eqno (4.7)
$$
having all the desirable features of (\ref{eq2.3}) and moreover the following one too,

\bigskip

\noindent (4.7.1) \quad For any $(x,y) \in M^2 (G_1 \mid Y(0)^{(0)})$ there is now a zipping path of length $\leq \Vert \Delta \Vert + 1$, in $Y(1)^{(1)}$. Assuming here that
$$
x \in H_n^0 (\gamma_1) \, , \ y \in H_n^0 (\gamma_2) \, ,
$$
then there will be a zipping path which makes use of the features (4.5) and (4.6), suggested in the diagram below
$$
\xymatrix{
&H_2^0(\gamma_1) \ar@{-}[rr]_{\lambda (h_2^0 , h_n^0 , \gamma_1)} &&H_n^0 (\gamma_1) \\
H_1^0(\gamma_0) \ar@{-}[dr]_{H_i^1 (\delta (\gamma_2))} \ar@{-}[ur]_{H_i^1 (\delta (\gamma_1))} \\
&H_2^0(\gamma_2) \ar@{-}[rr]_{\lambda (h_2^0 , h_n^0 , \gamma_2)} &&H_n^0 (\gamma_2).
}
$$
This {\ibf is} our zipping path of length $\leq \Vert \Delta \Vert + 1$ for $(x,y) \in M^2 (f)$.

\smallskip

A priori, the diagram of $0$-handles and continuous paths of $1$-handles (and $0$-handles) above, lives at the level of $\Delta$ but then, via equivariance it is transported all over the place. 

\smallskip

In the diagram above (part of the claim (4.7.1)) all the explicitly occuring $H^0$'s are $0$-handles of our original $Y(0)$. This also means that, so far, the zipping length is totally uncontrolled when we move to the additional $0$-handles which have been introduced by $Y(0) \Rightarrow Y(1)$.

\smallskip

But then, we can iterate indefinitely the $(Y(0) \Rightarrow Y(1))$-type construction, keeping all the time the same fundamental domain $\Delta \subset \tilde M^3 (\Gamma)$ as the first time. Moreover, like in \cite{26} and/or like in section III above, all the newly created $\{ \partial H$ and $\delta H \}$ will all the time be pushed closer and closer to the respective infinities. So, at least at each individual level $n$, local finiteness will stay all the time with us. Also, at each step $n$ we force equivariance and then $\{$restore GSC$\} + \{$decantorianize$\}$, a process which as already said keeps both the local finiteness and the equivariance alive. More explicitly, here is how our iteration will proceed. At each stage $n$ we will find some equivariant representation
$$
Y(0) \subset Y(1) \subset Y(2) \subset \ldots \subset Y(n-1) \subset Y(n)  \overset{G_n}{-\!\!\!-\!\!\!\longrightarrow} \tilde M^3 (\Gamma) \, , \eqno (4.8.n)
$$
such that for all $m < n$ we have $G_n \mid Y(m) = G_m$. The representation $Y(n)  \overset{G_n}{-\!\!\!-\!\!\!\longrightarrow} \tilde M^3 (\Gamma)$ has all the desirable features of (\ref{eq2.3}) from the lemmas~3.1, 3.2, 3.3 and, also in addition the following two too.

\bigskip

\noindent (4.9.$n$) \quad For any $(x,y) \in M^2 (G_n \mid Y(n-1)^{(0)})$ there is a zipping path of length $\leq \Vert \Delta \Vert + 1$, in $M^2 (G_n \mid Y(n)^{(1)})$.

\bigskip

\noindent (4.10) \quad Consider, like in the beginning of (4.9.$n$), an $(x,y) \in M^2 (G_n \mid Y(m-1)^{(0)})$ for some $m < n$. Then the zipping path which (4.9.$n$) assigns for this particular $(x,y)$, is exactly the one already assigned by (4.9.$m$).

\smallskip

With such compatibilities we have by now a well-defined space $Y(\infty)$ with a non-degenerate map
$$
Y(\infty) \underset{\rm def}{=} \bigcup_{n=0}^{\infty} Y(n) \overset{g(\infty) \underset{\rm def}{=} \lim G_n}{-\!\!\!-\!\!\!-\!\!\!-\!\!\!-\!\!\!-\!\!\!-\!\!\!-\!\!\!-\!\!\!-\!\!\!\longrightarrow} \tilde M^3 (\Gamma) \, , \eqno (4.10.1)
$$
and the issue now are the virtues of this (4.10.1).

\smallskip

To begin with, we have $Y(0) \in {\rm GSC}$. Next, our construction is such that the $1$-handles which via the inductive step $Y(n-1) \Rightarrow Y(n)$ are, to begin with, thrown in, are afterwards cancelled by later $2$-handles, thrown in by the same $Y(n-1) \Rightarrow Y(n)$. So $Y(\infty) \in {\rm GSC}$. The fact that, at each finite level $n$ we have $\Psi = \Phi$, easily implies this property for $g(\infty)$ too. So (4.10.1) is a representation, to begin with. At each finite $n$ we have equivariance, and zipping length is uniformly bounded by an $M$ which is independent of $n$. So (4.10.1) has these features too. But, what is, a priori, NOT automatic for $Y(\infty)$ is the local finiteness. This means that some care is required during the infinite process, so as to insure the convergence. But the technology from \cite{26} and/or the preceeding section saves here the day. Without any harm we can also add the following ingredient to our whole construction. 

\smallskip

Let us consider the generic bicollared handle $H_i^{\lambda} (\varepsilon , n)$ which is to be added to $Y(n-1)$ so as to get $Y(n)$. Let us also denote by $k(\lambda , i , \varepsilon , n)$ the level at which we attach it, in terms of the corresponding outgoing collar. Then, the whole infinite construction can be arranged so that we should have
$$
\lim_{n = \infty} \underset{\overbrace{\lambda , i , \varepsilon}}{\rm inf} k(\lambda , i , \varepsilon , n) = \infty \eqno (4.11)
$$
and also
$$
\lim_{\overbrace{p+i+\varepsilon + n = \infty}} \!\!\!\!\!\!\!\!\!\! (g(\infty) \, \delta H_{p,i}^{\lambda} (\varepsilon , n)) \subset g(\infty) \, \delta Y(\infty)^{(\lambda)} = (\tilde M^3 (\Gamma))^{(\lambda)} \, . \eqno (4.12)
$$
With this, our (4.10.1) has by now all the features demanded by lemma~4.1. So, by now the proof of the lemma in question is finished and then the proof of the $3^{\rm d}$ representation theorem too.

\smallskip

We will end this section with some remarks and comments.

\medskip

A) We very clearly have, in the context of our $3$-dimensional representation that $g(\infty) \, Y(\infty) = {\rm int} \, \tilde M^3 (\Gamma) = \tilde M^3 (\Gamma) - \partial \tilde M^3 (\Gamma)$, {\it i.e.} the image of $g(\infty)$ is locally finite. Moreover none of the special features like locally finite source, equivariance or bounded zipping length is involved here. But then, our $3^{\rm d}$ representation theorem is only a preliminary step towards the $2^{\rm d}$ representation theorem, stated in the introduction to this present paper and from which the proof of $\forall \, \Gamma \in {\rm GSC}$ will proceed, afterwards, in the subsequent papers of this series.

\smallskip

In the context of the $2$-dimensional representation (\ref{eq0.5}) we have again a locally finite source and also a uniformly bounded zipping length which, in some sense at least, should mean that ``all the zipping of (\ref{eq0.5}) can be performed in a finite time''. All these things notwithstanding and contrary to what one may a priori think, the image of the representation, {\it i.e.} the $fX^2$ FAILS NOW TO BE LOCALLY FINITE. The rest of this comment should give an idea of the kind of mechanism via which this seemingly self-contradictory phenomenon actually occurs.

\smallskip

In the same spirit as the diagram from (4.7.1) we consider inside $Y(\infty)$ and living over a same $\gamma \Delta \subset \tilde M^3 (\Gamma)$, $\gamma \in \Gamma$, infinitely many paths of $0$-handles and $1$-handles, all with the same $g(\infty)$-image, and starting at the same $H^0$
$$
H^0 \underset{\lambda_i}{\longrightarrow} H^0 (\gamma_i) \, ; \ i = 0,1,2,\ldots \eqno (A_1)
$$
We fix now our attention on $H^0 (\gamma_0)$ and consider some bicollared $1$-handle of $Y(\infty)$, call it $H^1$ attached to $H^0 (\gamma_0)$ in a direction transversal to the one of the the path $\lambda_0$ which hits $H^0 (\gamma_0)$. Some $2$-handle $H^2$ is now attached to a closed finite necklace $\ldots \cup H^0 (\gamma_0) \cup H^1 \cup \ldots$, inside $Y(\infty)$. We are here in a bicollared context meaning that $\partial H^2$ goes deep inside $H^0 (\gamma_0) \cup H^1$ and $g(\infty) (\delta H^0 (\gamma_0) \cup \delta H^1)$ cuts through $g(\infty) H^2$.

\smallskip

We move now to the $2$-dimensional representation (\ref{eq0.5}). Without going here into any more details, each $X^2 \mid H^{\lambda} (\gamma)$ is now a very dense $2$-skeleton of $H^{\lambda} (\gamma)$, becoming denser and denser as one approaches $\delta H^{\lambda} (\gamma)$, which is at infinity, as far as $H^{\lambda} (\gamma)$ is concerned.

\smallskip

When we go to $X^2 \mid H^2$, then among other things this will contain $2$-dimensional compact walls $W^2$ glued to $X^2 \mid (H^0 (\gamma_0) \cup H^1)$, parallel to the core of $H^2$, and such that we find a transversal intersection
$$
[f (\delta H^0 (\gamma_0) \cup \delta H^1) = g(\infty) (\delta H^0 (\gamma_0) \cup \delta H^1)] \pitchfork W^2 \subset \tilde M^3 (\Gamma) \, . \eqno (A_2)
$$
From $(A_1) + (A_2)$ we can pull out a sequence of double points $(x_n , y_n) \in M^2 (f)$, with $n=1,2,\ldots$ such that

\medskip

a) $x_n \in X^2 \mid H^0 (\gamma_n) \, , \ y_n \in W^2$.

\medskip

b) Inside $W^2$ we have $\lim y_n = y_{\infty} \in f (\delta H^0 (\gamma_0)) \cup W^2$, while as far as $g(\infty) H^0 (\gamma_n)$ is concerned, $\lim x_n = \infty$.

\medskip

c) Making use just of what $(A_1)$ and $W^2$ give us, we can produce zipping paths $\Lambda_n$ of length $\leq \Vert \Delta \Vert + 1$ connecting each $(x_n , y_n)$ to some singularity $\sigma_n \in X^2 \mid H^0$.

\medskip

d) Both $\sigma_n$ and $\Lambda_n$ go to infinity, at least in the following sense: inside $\tilde M^3 (\Gamma)$, the $f \Lambda_n$'s come closer and closer to $\underset{i}{\sum} \, f \partial H_i^0 \cup \underset{j}{\sum} \, f \partial H_j^1$. [This ``going to infinity'' of $\Lambda_n$ gets a more interesting meaning once one goes to the $S_u \, \tilde M^3 (\Gamma)$ from the introduction and one applies punctures, like in (\ref{eq0.9}).]

\medskip

But the point here is that $fX^2$ is {\ibf not locally finite} in the neighbourhood of $fy_{\infty}$.

\medskip

B) Notice also, that contrary to what one might have thought, a priori, bounded zipping length and difficulty of solving algorithmic problems for $\Gamma$, like the word problem (or lack of such difficulties), are totally unrelated to each other.

\medskip

C) There is no finite uniform boundedness for the lengths of the zipping paths induced by (4.8.9) on $Y (n)^{(0)} - Y(n-1)^{(0)}$. That is why an INFINITE CONSTRUCTION is necessary here.

\medskip

D) So, eventually, our uniform boundedness has been achieved via an {\ibf infinite} process which is {\ibf convergent}. Also, it would not cost us much to put also some lower bound for the zipping lengths too, if that would be necessary.

\medskip

E) The way our proof proceeded, was to show (4.2) $+$ (4.3) for any $(x,y) \in M^2 (G_{\infty}) \mid \{\mbox{$0$-skeleton of} \ Y(\infty)\}$. But from then on, it is very easy to connect, first, any $(x,y) \in M^2 \mid \{\mbox{$1$-skeleton}\}$ to a neighbouring double point in the $0$-skeleton, via a short zipping path confined in the $1$-skeleton. Next any double point involving the $2$-skeleton can be connected, via short zipping paths too, to neighbouring double points in the $0$-skeleton.

\medskip

F) Our general way of proceeding, in this paper, was the following. We have started with the initial representation $X \overset{f}{\longrightarrow} \tilde M^3 (\Gamma)$ from section I (see for instance (2.9.1)), then we enlarged it into the $Y \overset{G}{\longrightarrow} \tilde M^3 (G)$ from (\ref{eq2.3}) and then finally to our present gigantic $Y(\infty) \overset{g(\infty)}{-\!\!\!-\!\!\!-\!\!\!-\!\!\!\longrightarrow} \tilde M^3 (\Gamma)$ (\ref{eq3.1}). We started with the initial feature $X \in {\rm GSC}$ and thus we never lost it in our successive constructions.

\medskip

G) When, in the next paper, we will change our (4.1) into a 2-dimensional representation, then the (4.12) will make that the second finiteness condition is again satisfied.

\end{document}